\documentclass[journal]{IEEEtran}
\pdfoutput=1       
\usepackage[cmex10]{amsmath}
\usepackage{algorithm}
\usepackage{algorithmic}
\usepackage{array}
\usepackage{mdwmath}
\usepackage{mdwtab}
\usepackage{eqparbox}
\usepackage{float}
 \usepackage[font=footnotesize]{subfig}
\usepackage{fixltx2e}
\usepackage{stfloats}
\usepackage{amsfonts}
\usepackage{bbm}
\usepackage{amssymb}
\usepackage{cite}
\usepackage{graphicx}

\usepackage[usenames]{color}

\newtheorem{theorem}{\bf Theorem}%[section]
\newtheorem{lemma}{Lemma}%[section]
\newtheorem{corollary}{Corollary}%[section]
\def\ba{\begin{array}}
\def\ea{\end{array}}
\newcommand{\beq}{\begin{equation}}
\newcommand{\eeq}{\end{equation}}
\newcommand{\bq}{\begin{eqnarray}}
\newcommand{\eq}{\end{eqnarray}}
\newcommand{\bqn}{\begin{eqnarray*}}
\newcommand{\eqn}{\end{eqnarray*}}
\newcommand{\bee}{\begin{enumerate}}
\newcommand{\eee}{\end{enumerate}}
\newcommand{\bi}{\begin{itemize}}
\newcommand{\ei}{\end{itemize}}

\newcommand{\mathN}{\mathcal{N}}
\newcommand{\mathE}{\mathcal{E}}
\newcommand{\mathG}{\mathcal{G}}

\newcommand{\droped}[1]{{\color{blue} \sout{}}}
\newboolean{showcomments}
\setboolean{showcomments}{false}
\newcommand{\slow}[1]{\ifthenelse{\boolean{showcomments}}
{ \textcolor{red}{(Steven says:  #1)}}{}}
\newcommand{\qiuyu}[1]{  \ifthenelse{\boolean{showcomments}}
{ \textcolor{blue}{(Qiuyu says:  #1)}} {}  }

\begin{document}

%\title{\LARGE \bf Optimal Branch Exchange for Feeder Reconfiguration in Distribution Networks}
\title{Optimal Branch Exchange for Distribution System Reconfiguration}

\author{Qiuyu Peng and Steven H. Low% <-this % stops a space 
%\\Engr. \& App. Sci., Caltech, CA
\thanks{*This work was supported by NSF NetSE grant CNS 091104, ARPA-E grant DE-AR0000226, Southern California Edison, National Science Council of Taiwan, R.O.C. grant NSC 103-3113-P-008-001 and Resnick Institute.}% <-this % stops a space
\thanks{Qiuyu Peng is with the Electrical Engineering Department 
and Steven H. Low is with the Computing and Mathematical Sciences and
the Electrical Engineering Departments, California Institute of Technology, Pasadena, CA 91125, USA.
{\small \tt \{qpeng, slow\}@caltech.edu}}%
}

\maketitle
\thispagestyle{empty}
\pagestyle{empty}

\begin{abstract}
The  feeder reconfiguration problem chooses the on/off status of the switches 
in a distribution network in order to minimize a certain cost such as power loss.   It is
a mixed integer nonlinear program and hence hard to solve.   A popular heuristic search
 consists of repeated application of branch exchange, where some loads are transferred 
from one feeder to another feeder while maintaining the radial structure of the network,
until no load transfer can further reduce the cost.
Optimizing each branch exchange step is itself a mixed integer nonlinear program.
In this paper we propose an efficient algorithm for optimizing a branch exchange step.
It uses an AC power flow model and is based on the recently developed convex relaxation 
of optimal power flow.
We provide a bound on the gap between the optimal cost and that of our solution.
We prove that our algorithm is optimal when the voltage magnitudes are the same at all buses.
We illustrate the effectiveness of our algorithm through the simulation of real-world distribution
feeders.
\end{abstract}

%------------------------------------------------------------------------------------------------------------------------------------------------------------------------
\begin{IEEEkeywords}
Distribution System, Feeder Reconfiguration, Nonconvex Optimization
\end{IEEEkeywords}

 \section{Introduction}

\subsection{Motivation}

\IEEEPARstart{P}{rimary} distribution systems  have sectionalizing switches that connect line sections
and tie switches that connect two primary feeders, two substation buses, or loop-type laterals.
In normal operation these switches are configured such that a distribution network is acyclic 
and every load is connected to exactly one substation.
The topology of the network can be reconfigured  by changing the on/off status of these
switches, for the purpose of load balancing, loss minimization, or service restoration, e.g., 
\cite{Merlin1975, Civanlar1988, Baran1989c, Chiang1990a}.  See also an survey in \cite{sarfi1994}
for many early papers and references to some recent work in \cite{Jabr2011}.

For instance when a single feeder is overloaded, a currently open tie switch can be closed to 
connect the feeder to another substation.  Since this will create a loop or connect some loads to
two substations, a currently closed sectionalizing switch will be opened to maintain a radial 
topology in which every load is connected to a single substation.
Following \cite{Baran1989c} we call this a ``branch exchange''  
where the goal is to select the pair of switches for closing/opening that 
achieves the best load balancing. 
More generally one can optimize a certain objective over the topology of the entire distribution
network by choosing the on/off status of all the switches, effectively selecting a best spanning 
tree among all possible spanning trees of the network topology.  
Even though the problem of minimum spanning tree has been well studied \cite{Cormen2001},
the problem here is different.  Unlike the standard minimum spanning tree problem where the link
costs are fixed and the minimization is only over the topology, in our case, the link costs result
from an optimal power flow (OPF) problem that must be solved for each candidate spanning tree.
This is therefore a mixed integer nonlinear program and can generally be NP-hard. 
As a result the large majority of proposed solutions are heuristic in nature \cite{sarfi1994}; see
also references in \cite{Jabr2011}.
A heuristic search method is proposed in \cite{Civanlar1988, Baran1989c} which we discuss in more detail
below.
The problem is formulated as a multi-objective mixed integer constrained optimization in \cite{Chiang1990a}
and solved using a simulated annealing technique.
Ordinal optimization is proposed in \cite{Jabr2011} to reduce the computational burden through
order comparison and goal softening.
Unlike these heuristic methods, an interesting exhaustive search method is proposed in \cite{Morton2000} to compute a
globally optimal solution under the assumption that loads are constant-current, instead of constant-power
as often assumed in load flow analysis.  Starting from an initial spanning tree,
the proposed method applies the branch exchange step % of \cite{Civanlar1988}
in a clever way to generate all spanning trees exactly once and efficiently compute the power loss for
each tree recursively in order to find a tree with the minimum loss.
A constant-current load model is also used in \cite{Liu1989} where the optimization problem becomes a
mixed integer {\em linear} program.   A global optimality condition is
derived and an algorithm is provided that attains global optimality under certain conditions.
Recently sparse recovery techniques have been applied to this problem in 
\cite{Dall'Anese2013-reconfig1, Dall'Anese2013-reconfig2} where the network is assumed to be unbalanced
and the optimization is formulated in terms of currents.  
  Even though the Kirchhoff current law at each node is linear, nonconvexity arises due to the 
binary variables that represent the status of the switches and the quadratic (bilinear) relation between the
power injection (representing load or generation) and the current injection at each bus.
To deal with the former nonconvexity,
\cite{Dall'Anese2013-reconfig1} removes the binary variable and adds a regularization term that encourages
\emph{group sparsity}, i.e., at optimality, either a branch current is zero for all phases (corresponding to opening
the line switch) or nonzero for all phases (closing the line switch).
The latter nonconvexity is removed by approximating the quadratic relation by a linear relation.  The resulting
approximate problem is a (convex) second-order cone problem and can be solve efficiently.
The formulation in \cite{Dall'Anese2013-reconfig2} adds a chance constraint that the probability of loss of load is
less than a threshold.  Chance constraints are generally intractable and the paper proposes to replace it
with scenario-based approximation which is convex.

In this paper we study a single branch exchange step proposed in \cite{Civanlar1988}.
Each step transfers some loads from one feeder to another if it reduces the overall cost.
An efficient solution for a single branch exchange step is important because, as  suggested in \cite{Baran1989c},
a heuristic approach for optimal network reconfiguration consists of repeated application of branch exchanges
until no load transfer between two feeders can further decrease the cost.
This simple greedy algorithm yields a local optimal.
The key challenge is to estimate the cost reduction for each load transfer. 
Specifically once a currently open tie switch has been selected for closing, the issue is to determine
which one of several currently closed sectionalizing switches should be opened that will provide the largest 
cost reduction.
Each candidate sectionalizing switch (together with the given tie switch) transform the existing spanning
tree into a new spanning tree.
A naive approach will solve an OPF for each of the candidate spanning tree and choose one
that has the smallest cost.
This may be prohibitive both because the number of candidate spanning trees can be large
and because OPF is itself a nonconvex problem and therefore hard to solve.
The focus of \cite{Civanlar1988, Baran1989c, Wang1995} is to develop much more efficient ways
to approximately evaluate the cost reduction by each candidate tree without solving the full power flow problem.
The objective of \cite{Civanlar1988} is to minimize loss and it derives a closed-form expression for
approximate loss reduction of a candidate tree.  This avoids load flow calculation altogether.
A new branch flow model for distribution systems is introduced in \cite{Baran1989c} that allows a
recursive computation of cost reduction by a candidate tree.
This model is extended to unbalanced systems in \cite{Wang1995}.

\subsection{Summary}

We make two contributions to the solution of branch exchange.
First we propose a new algorithm to determine the sectionalizing switch whose opening
will yield the largest cost reduction, once a tie switch has been selected for closing.  
We use the full AC power flow model introduced in \cite{branch1,branch2} for radial system and solve them through the method of convex relaxation developed recently in \cite{Farivar-2011-VAR-SGC, Farivar-2013-BFM-TPS, Gan-2013-BFMt-CDC}.
Instead of assuming constant real/reactive power for each load bus as prior works, we consider the scenario where the real/reactive power can also be control variables.
Moreover the algorithm requires solving at most three OPF problems regardless of the number
of candidate spanning trees.
Second we bound the gap between the cost of our algorithm and the optimal cost.
We prove that when the voltage magnitude of each bus is the same our algorithm is optimal.
We illustrate our algorithm on two Southern California Edition (SCE) distribution feeders, and
in both cases, our algorithm has found the optimal branch exchange.

The rest of the paper is organized as follows. We formulate the optimal feeder reconfiguration 
problem in Section \ref{sec:model} and propose an algorithm to solve it in Section \ref{sec:alg}. 
The performance of the algorithm is analyzed  in Section  \ref{sec:analysis}. 
The simulation results on SCE distribution circuits are given in Section \ref{sec:simulation}. 
We conclude in Section \ref{sec:conclusion}.

%----------------------------------------------------------------------------------------------------------------
\section{Model and problem formulation}\label{sec:model}

 \begin{figure}
\centering
\includegraphics[width=0.9\linewidth]{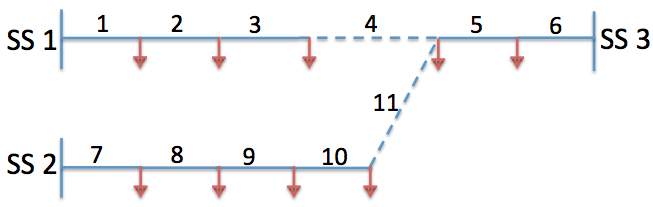}
\caption{A distribution network. 
Solid Lines are closed and dash lines are open. The red arrows are load buses.}
\label{fig:feeder1}
\end{figure}

\subsection{Feeder reconfiguration}
\label{sec:background}

A distribution system consists of buses, distribution lines, and (sectionalizing and tie)
switches that can be opened or closed.   There are two types of buses.  {\em Substation
buses} (or just {\em substations}) are connected to a transmission network from which 
they receive bulk power, and {\em load buses} that receive power from the substation
buses.   During normal operation the switches are configured so that 
\bee
\item The network is radial, i.e., has a tree topology.
\item Each bus is connected to a single substation.
\eee
We will refer to the subtree that is rooted at each substation
bus as a feeder; hence each feeder is served by a single substation.
Optimal feeder reconfiguration is the problem of reconfiguring the switches to minimize
a certain cost subject to the two constraints above, in addition to operational constraints
on voltage magnitudes, power injections, and line limits.

We assume that there is an on/off switch on each line (i.e., modeling the subsystem between
each pair of switches as a single line),
 and focus on an iterative greedy algorithm first proposed in \cite{Baran1989c}. 
 We illustrate this algorithm on the simple network shown in Fig. \ref{fig:feeder1}
 where solid and dash lines represent closed and open switches respectively. 

 There are $3$ feeders, each of which connects to one substation, SS1, SS2, or SS3. 
 Suppose lines 4 and 11 are open in the current iteration.  In each iteration one
 of the open switches is selected and closed, say, that on line 4.  
 This joins two feeders so that every bus along lines 1 to 6 are now connected
 to both substations SS1 and SS3.    To restore the property that each bus is connected
 to a single substation,
 we then choose one line among $\{1,2,3,4,5,6\}$ to open that minimizes the cost.
 This two-step procedure is called a branch exchange.
This procedure is repeated until the configuration stabilizes, i.e., the line that is chosen
to open in step two is the original open line selected in step one, for all open switches.
In summary, each iteration of the algorithm consists of two steps:
\bee
\item  Chooses a line $e_1$ with an open switch and close the switch. 
\item  Identify a line $e_2$ in the two feeders that was joined in Step 1
	to open that minimizes the objective. 
\eee
The algorithm terminates when $e_1=e_2$ for all the open switches. 
The greedy search only guarantees a local optimum since it may terminate before
searching through all spanning trees. 
In this paper we propose an efficient and accurate method to accomplish Step 2
in each branch exchange (iteration).
We will use the nonlinear (AC) power flow model and apply convex relaxations
developed recently for its solution.
Most existing algorithms that we are aware of perform Step 2 based on either linearized 
power flow equations \cite{Baran1989c,Civanlar1988} or assumption of constant current loads \cite{Morton2000}.   Linearized power flow (called DC power flow) model is reasonable in transmission networks but is less so in distribution networks.

%\slow{Is the following true?   I haven't read many papers, so I'm not sure.  
%Linearized DistFlow is used in \cite{Baran1989c}, though that's only one of 2 methods
%proposed there.  Is the expression in \cite{Civanlar1988} also a sort of linearization?
%\cite{Morton2000} has a linear model but that's not an approximation, but results from
%its assumption of constant-current loads.  }
%\qiuyu{\cite{Civanlar1988} is another sort of Linearization, it can be seen as kind of DC-PF. When they transfer loads, their estimate of the $\Delta P$ is a linearized estimate based on the current physical variable before the load transfer.}

\begin{figure}
\centering
\subfloat[Two feeders served by different substations.] {
	\includegraphics[scale=0.25]{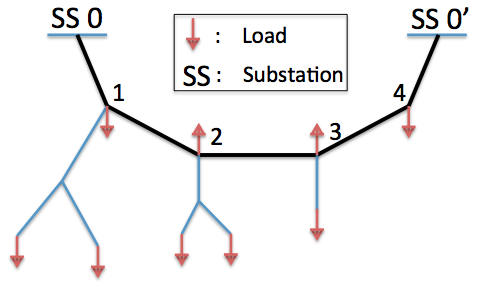}
	\label{fig:step1a}
	}
\hspace{0.2in}
\subfloat[Two feeders served by the same substation.]{
	\includegraphics[scale=0.25]{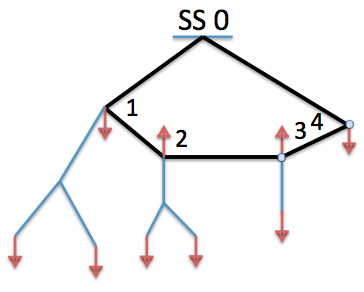}
	}
% \caption{Two newly connected feeders served by different substations (left panel)
%	or the same substation (right panel).}
\caption{Feeders after step 1 of a branch exchange.}
\label{fig:step1}
\end{figure}

After we close the switch on a line, there are two possible cases (see Fig \ref{fig:step1}): 
(1) The two connected feeders are served by different substations; or
(2) The two connected feeders are served by the same substation. 

In both cases the switch on one of the lines needs to be opened.   Case (2) can be reduced
to case (1) by replacing the substation 0 by two virtual substations $0$ and $0'$ 
as shown in  Fig. \ref{fig:step1a}. 

We now describe our model and  formulate the problem of determining the 
optimal switch to open along the path that connects two substations.

\subsection{Network model}\label{subsec:model}

\begin{figure}
\centering
\includegraphics[width=0.6\linewidth]{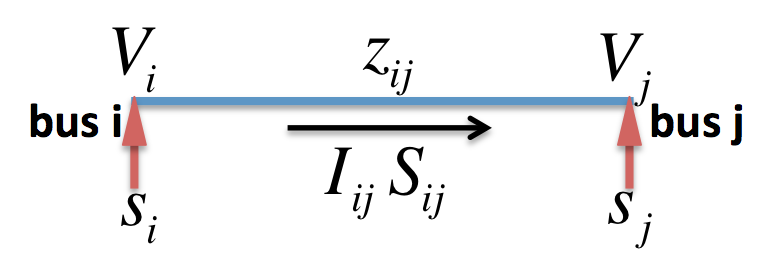}
\caption{Notations.}
\label{fig:model}
\end{figure}

We consider an AC power flow model where all variables are complex. 
A distribution network is denoted by a graph $\mathG(\mathN,\mathE)$, where nodes in $\mathN$ represent buses and edges in $\mathE$ represent distribution lines. For each bus $i\in \mathN$, let $V_i=|V_i|e^{\mathbf{i}\theta_i}$ be its complex voltage and $v_i:=|V_i|^2$ be its magnitude squared. 
Let $s_i=p_i+\mathbf{i}q_i$ be its net power injection which is defined as generation minus consumption. 
We associate a direction with each line $(i,j)\in \mathE$  represented by 
an ordered pair of nodes in $\mathN$. For each line $(i,j)\in \mathE$,
 let $z_{ij}=r_{ij}+\mathbf{i}x_{ij}$ be its complex impedance and $y_{ij}:=1/z_{ij}$  its admittance. We have $x_{ij}>0$ since lines are inductive. Let $I_{ij}$ be the complex branch current from 
 buses $i$ to $j$ 
and $\ell_{ij}:=|I_{ij}|^2$ be its magnitude squared. 
Let $S_{ij}=P_{ij}+\mathbf{i}Q_{ij}$ be the branch power flow from buses $i$ to $j$. 
For each line $(i,j)\in \mathE$, define $S_{ji}$ in terms of $S_{ij}$ and $I_{ij}$ 
by $S_{ji}:=-S_{ij}+\ell_{ij}z_{ij}$. Hence $-S_{ji}$ represents the power received 
by bus $j$ from bus $i$. 
The notations are illustrated in Fig. \ref{fig:model}. A variable without a subscript 
denotes a column vector with appropriate components, as summarized below.
\vspace{0.05in}

\begin{tabular}{| l | l | }
\hline
  $p:=(p_{i},i\in \mathN) $  &$q:=(q_{i},i\in \mathN) $  \\
  \hline
  $P:=(P_{ij},(i,j)\in \mathE)$   &$Q:=(Q_{ij},(i,j)\in \mathE)$  \\
  \hline
  $v:=(v_{i},i\in \mathN)$ & $\ell:=(\ell_{ij},(i,j)\in \mathE)$\\
  \hline
\end{tabular}
\vspace{0.05in}

For most parts of our paper (except the proof of Lemma \ref{lem:branchincrease}), it
suffices to work with a `relaxed' model, first proposed in \cite{branch1,branch2} to 
model radial network,  where we ignore the phase angles of voltages
and currents and use only the variables $x := (p, q, P, Q, \ell, v)$.   These variables 
satisfy:
\begin{align}
p_i&=-\sum_{(k,i)\in \mathE}\left(P_{ki}-\ell_{ki}r_{ki}\right)+\sum_{(i,j)\in \mathE}P_{ij}, \ \  i\in \mathN
\label{eqn:opf3}
\\
q_i&=-\sum_{(k,i)\in \mathE}\left(Q_{ki}-\ell_{ki}x_{ki}\right)+\sum_{(i,j)\in \mathE}Q_{ij}, \ \  i\in \mathN 
\label{eqn:opf4}
\\
v_j&=v_i -2(r_{ij}P_{ij}+x_{ij}Q_{ij})+\ell_{ij}|z_{ij}|^2, 	 \  (i,j) \in \mathE
\label{eqn:opf2}
\\
\ell_{ij} v_i &= {P^2_{ij}+Q^2_{ij}},  \ \  (i,j) \in \mathE
\label{eqn:lastopf}
\end{align}
Given a vector $x$ that satisfies \eqref{eqn:opf3}--\eqref{eqn:lastopf},
the phase angles of the voltages and currents can be uniquely determined 
for a radial network, and therefore this relaxed model \eqref{eqn:opf3}--\eqref{eqn:lastopf}
is equivalent to the full AC power flow model for a radial network;
See \cite[section III-A]{Farivar-2013-BFM-TPS} for details. 

%Ohm's Law for $(i,j)\in E$:
%\bq \label{eqn:ohm}
%V_i-V_j=I_{ij}z_{ij}
%\eq
%Definition of branch power for $(i,j)\in E$:
%\bq\label{eqn:branchpower}
%S_{ij}:=V_iI_{ij}^*
%\eq
%Power balance at each bus $i\in N$:
%\bq\label{eqn:nodebalance}
%s_i=-\sum_{(k,i)\in E}\left(S_{ki}-|I_{ki}|^2z_{ki}\right)+\sum_{(i,j)\in E}S_{ij}
%\eq
In addition $x$ must also satisfy the following operational constraints:
\bi
\item {\it Power injection constraints:} for each bus $i\in \mathN$
\bq\label{eqn:injection}
\underline{p}_i\leq p_i\leq \overline{p}_i \quad\text{and}\quad \underline{q}_i\leq q_i\leq \overline{q}_i
\eq

\item {\it Voltage magnitude constraints:} for each bus $i\in \mathN$
\bq\label{eqn:Branch_voltage}
\underline{v}_i \leq v_i \leq \overline{v}_i
\eq

\item {\it Line flow constraints:} for each line $(i,j)\in \mathE$
\bq\label{eqn:Branch_lineloss}
\ell_{ij} \leq {\overline \ell_{ij}}
\eq
\ei
\subsection{Problem formulation}

As described in section \ref{sec:background}, we focus on reconfiguring a
network path where two feeders are served by two different substations
as shown in Fig. \ref{fig:step1a}. 
Consider a (connected) tree network $\mathG(\mathN, \mathE)$.
 $\mathN:=\{0,1,\ldots,n,0'\}$ denote the set of buses, where the two substations 
are indexed by $0, 0'$ and the load buses are indexed by $\{1,\ldots,n\}$.

Since $\mathG$ is  a tree there is a unique path between any two 
buses in $\mathG$. For every pair of buses $i,j\in \mathN$, let 
$\mathE(i,j)\subseteq \mathE$ be the collection of edges on the unique path
between $i$ and $j$.
Given any subgraph $\mathG'$ of $\mathG$ let $x^{\mathG'} := (p^{\mathG'}, q^{\mathG'}, P^{\mathG'}, Q^{\mathG'},
\ell^{\mathG'}, v^{\mathG'})$ denote the set of variables defined on $\mathG'$ with appropriate
dimensions.  For notational
simplicity we often ignore the superscript $\mathG'$ and write $x := (p, q, P, Q, \ell, v)$
instead when the meaning is clear from the context. 
Given any subgraph $\mathG'$ of $\mathG$, let 
$\mathbb X(\mathG') :=\{x^{\mathG'}\mid x^{\mathG'} \mbox{ satisfies } \eqref{eqn:opf3}-\eqref{eqn:Branch_lineloss}\}$
be the feasible set of variables $x$ defined on $\mathG'$. 
In particular, $\mathbb X(\mathG)$ is the feasible set for the entire distribution network 
represented by $\mathG$.

Each load bus connects to one substation if and only if one of the switch on a line in path $\mathcal{E}(0,0')$ is turned off.
Given a (connected) tree $\mathG(\mathN,\mathE)$ and a path $\mathcal{E}(0,0')$ between bus 
$0$ and $0'$, denote by $\mathG_0^i(\mathN_0,\mathE_0)$ and $\mathG_{0'}^j(\mathN_{0'},\mathE_{0'})$ 
 the two subtrees after we remove line $(i,j)\in\mathcal{E}(0,0')$, where $0\in \mathN_0$ and 
 $0'\in \mathN_{0'}$. The  minimum power injections for $\mathG_0^i$ and $\mathG_{0'}^j$ are 
 defined as 
\bq
p_0^{i}&:=&\min_{x\in\mathbb{X}(\mathG_0^i)} p_0 \label{eqn:opfi}\\
p_{0'}^{j}&:=&\min_{x\in\mathbb{X}(\mathG_{0'}^j)} p_{0'} \label{eqn:opfj}
\eq

The optimal branch exchange for feeder reconfiguration problem is defined as:
\bqn
\text{\textbf{OFR-branch} (OFR): } \min_{(i,j)\in \mathcal{E}(0,0')} \Gamma(p_0^{i},p_{0'}^{j})
\eqn
where $\Gamma(p_0,p_{0'})$ can be any convex increasing cost function. 
When $\Gamma(p_0,p_{0'})=p_0+p_{0'}$, our goal is to minimize the aggregate power 
injection from the substations.   Since $p_0+p_0'$ equals the aggregate load
(real power consumption) in the network and the total real power loss, if the loads
are fixed, then minimizing $p_0 + p_0'$ also minimizes power loss.
For simplicity we will also refer to OFR-branch as OFR in this paper.

A naive solution to  OFR is to enumerate all the lines in 
$\mathcal{E}(0,0')$ and compare the objective value for each case.    It 
is inefficient as it requires solving two optimal power flow (OPF) problems 
\eqref{eqn:opfi} and \eqref{eqn:opfj} for each line.   This can be computationally
expensive if the size of $\mathcal{E}(0,0')$ is large.
In the following we will develop an algorithm to solve OFR that involves solving
at most three OPF problems regardless of the size of $\mathcal{E}(0,0')$. 

The OPF problem is itself a nonconvex problem thus even one OPF problem is hard to solve in general. However, the OPF problems involved in 
the proposed algorithm can be solved through a convex relaxation. Next, we briefly describe SOCP (second-order cone program) relaxation of 
OPF recently developed in \cite{Farivar-2011-VAR-SGC, Farivar-2013-BFM-TPS, Gan-2013-BFMt-CDC,V_opf}.

\subsection{OPF and convex relaxation}\label{subsec:opf}

The optimal power flow problem seeks to optimize a certain objective 
over the feasible set $\mathbb{X}(\mathG)$ specified by the power flow equations 
\eqref{eqn:opf3}-\eqref{eqn:lastopf} and the operation constraints \eqref{eqn:injection}-\eqref{eqn:Branch_lineloss}:
\bqn
\text{\textbf{OPF-$\mathG$}:}\quad\min_{x\in  \mathbb{X}(\mathG)}\Gamma(p_0,p_{0'})
\eqn
It is a non-convex problem due to the quadratic equalities \eqref{eqn:lastopf}. 
Relaxing \eqref{eqn:lastopf} to inequalities:
\bq
\ell_{ij} v_i  &\geq& {P^2_{ij}+Q^2_{ij}}
\label{eqn:lastopf1}
\eq
leads to a second order cone program (SOCP) relaxation. 
Formally define  $\mathbb{X}_c(\mathG):=\{x\mid x \mbox{ satisfies } \eqref{eqn:opf3}- \eqref{eqn:opf2},\eqref{eqn:injection}-\eqref{eqn:Branch_lineloss},\eqref{eqn:lastopf1}\}$. 
The SOCP relaxation of OPF-$\mathG$ is:  
\bqn
\text{\textbf{SOPF-$\mathG$}:}\quad\min_{x\in \mathbb{X}_c(\mathG)}\Gamma(p_0,p_{0'})
\eqn 
SOPF-$\mathG$ is convex and can be solved efficiently.
Clearly SOPF-$\mathG$ provides a lower bound for OPF-$\mathG$ since 
$\mathbb{X}\subseteq \mathbb{X}_c$.
It is called {\it exact} if every solution $x^*$ of SOPF-$\mathG$ attains  equality in 
\eqref{eqn:lastopf1}. 
For radial  networks  SOCP relaxation is exact under some mild conditions 
\cite{Farivar-2011-VAR-SGC, Farivar-2013-BFM-TPS, Gan-2013-BFMt-CDC,V_opf}. 

Since the network graph is radial after we join two feeders severed by different substations as shown in Fig. \ref{fig:step1a}, we will assume that the SOCP relaxation of OPF is always exact throughout this paper.
In that case we can solve 
SOPF-$\mathG$ and recover an optimal solution to the original non-convex OPF-$\mathG$. 
A similar approach can be applied to the OPF problems defined in \eqref{eqn:opfi} and \eqref{eqn:opfj}.

%------------------------------------------------------------------------------------------------------------------------------------------------------------------------
\section{Algorithm of Branch exchange for Feeder Reconfiguration}\label{sec:alg}

 OFR  seeks to minimize $\Gamma(p_0,p_{0'})$ by opening the switch on a line in $\mathcal{E}(0,0')$. 
 Let $k_0,k_{0'}\in \mathN$ denote the buses such that $(0,k_0),(k_{0'},0')\in \mathcal{E}(0,0')$. The algorithm for OFR is given in Algorithm 1.

\begin{table}[htbp]
\centering
\normalsize
\begin{tabular}{|p{8.2cm}|}
\hline
{\bf Algorithm 1}: Branch Exchange Algorithm for OFR\\
\hline
{\bf Input}: objective $\Gamma(p_0,p_{0'})$,
network constraints $(\overline p,\underline p, \overline q,\underline q,\overline \ell,\overline v,\underline v)$.\\
{\bf Output}: line $e^*$.\\
\hline
Solve OPF-$\mathG$; let $x^*$ be an optimal solution.\\
\bee
\item $P^*_{0,k_0}\leq 0$: $e^*\leftarrow(0,k_0)$.
\item $P^*_{k_{0'},0'}\geq 0$: $e^*\leftarrow(k_{0'},0')$.
\item $\exists (k_1,k_2)\in \mathcal{E}(0,0')$ such that $P^*_{k_1k_2}\geq 0$ and $P^*_{k_2k_1}\geq 0$: $e^*\leftarrow(k_1,k_2)$.
\item $\exists (k_1,k_2), (k_2,k_3)\in \mathcal{E}(0,0')$ such that $P^*_{k_2,k_1}\leq 0$  and $P^*_{k_2,k_3}\leq 0$. Calculate $p_0^{k_1}$, $p_{0'}^{k_1}$, $p_0^{k_2}$ and $p_{0'}^{k_3}$.
     \begin{itemize}
\item $\Gamma(p_0^{k_1},p_{0'}^{k_2})\geq\Gamma(p_0^{k_2},p_{0'}^{k_3})$: $e^*\leftarrow(k_2,k_3)$.
\item $\Gamma(p_0^{k_1},p_{0'}^{k_2})<\Gamma(p_0^{k_2},p_{0'}^{k_3})$: $e^*\leftarrow(k_1,k_2)$.
\end{itemize}
\eee\\
\hline
\end{tabular}
\end{table}

 The basic idea of Algorithm 1 is simple and we illustrate it using the line 
 network in Fig. \ref{fig:linenetwork}. After we solve OPF-$\mathG$ with $x^*$:
\bee
\item if bus $0$ receives positive real power from bus $1$ through line $(0,1)$, open line  $(0,1)$.
\item if bus $0'$ receives positive real power from bus $n$ through line $(n,0')$, open line  $(n,0')$.
\item if there exists a line $(k,k+1)$ where positive real power is injected from both ends,
	open line $(k,k+1)$.
\item if there exists a bus $k$ that receives positive real power from both sides, open either
	line $(k-1,k)$ or $(k,k+1)$.
\eee 

We are interested in the performance of Algorithm 1, specifically:
\begin{itemize}
\item Is the solution $x^*$ to OPF-$\mathG$ unique and satisfies exactly one of the cases
 $1)-4)$? 
\item Is the line $e^*$ returned by Algorithm 1 optimal for OFR?
\end{itemize}
We next state our assumptions and answer these two questions under those assumptions.

\section{Performance of Algorithm 1}\label{sec:analysis}

\begin{figure}
\centering
\includegraphics[width=0.9\linewidth]{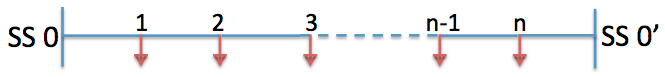}
\caption{A line Network}
\label{fig:linenetwork}
\end{figure}

For ease of presentation we only prove the results for a line network as shown
in Fig. \ref{fig:linenetwork}. They generalize in a straightforward manner to
radial networks as shown in Fig. \ref{fig:step1a}.

Our analysis is divided into two parts. First we show that, OPF-$\mathG$ has a 
unique solution $x^*$ and it satisfies exactly one of the cases 1) - 4) in Algorithm 1.   
This means that Algorithm 1 terminates correctly.   Then we
 prove that the performance gap between the solution $e^*$ given by Algorithm 1 
and an optimum of OFR is zero when the voltage magnitude of every bus is fixed at the
same nominal value, and bound the gap  by a small value when the voltage
magnitudes are fixed but different.

\subsection{Assumptions}

For the line network  in Fig. \ref{fig:linenetwork}, let the buses at the two ends
be substation buses and buses in between be load buses.
Hence the path between substations $0$ and $0'$ is $\mathcal{E}(0,0')=\mathE$;
we sometimes use $0'$ and $n+1$ interchangeably for ease of notation.
We collect the  assumptions we need as follows:
\bee
\item [A1]: $\overline{p}_{k}<0$ for $1\leq k\leq n$ and $\overline{p}_{k}>0$ for $k=0,0'$.
\item [A2]: $\overline v_k=\underline v_k$, $\overline q_k=-\underline q_k=\infty$ for $k\in \mathN$.
\item [A3]: $|\theta_i-\theta_j| <\arctan (x_{ij}/r_{ij})$ for $(i,j)\in \mathE$.
\item [A4]: The feasible set $\mathbb{X}(\mathG)$ is compact.
\eee
A1 is a key assumption and it says that
buses 0 and $0'$ are substation buses that inject positive real power
while buses $1, \dots, n$ are load buses that absorb real power.
A2 says that the voltage magnitude at each bus is fixed at their nominal value. 
To achieve this we also require that the reactive power injections are unconstrained. 
This is a reasonable approximation for our purpose
since there are Volt/VAR control mechanisms on distribution networks that
maintain voltage magnitudes within a tight range around their nominal values
as demand and supply fluctuate \cite{Farivar-2011-VAR-SGC}. Our simulation results on real SCE feeders show 
that the algorithm also works well without A2.
 A3 is a technical assumption that bounds the angle difference between adjacent buses. It, together with A2, guarantees that SOCP relaxation of OPF is exact \cite{V_opf}.\footnote{Although voltage phase angles $\theta_i$ are
 relaxed in the relaxed branch flow model \eqref{eqn:opf3}-\eqref{eqn:lastopf},
 they are uniquely determined by $\theta_i-\theta_j=\angle (v_i - z_{ij}^* S_{ij})$ in
 a radial network \cite{Farivar-2013-BFM-TPS}.}
A4 is an assumption that is satisfied in practice and guarantees
that our optimization problems are feasible.

\subsection{Main results}

Algorithm 1 needs to solve up to three OPF problems.
The result of \cite{V_opf}  implies that we can solve these problems
through their SOCP relaxation.
\begin{theorem}\label{thm:exact}
Suppose A2 and A3 hold.    Then, for any subgraph $\mathG'$ of $\mathG$ (including
$\mathG$ itself),
\bee
\item SOPF-$\mathG'$ is exact provided the objective function $\Gamma(p)$ is a
convex nondecreasing function of $p$.
\item OPF-$\mathG'$ has a unique solution provided the objective function $\Gamma(p)$ is convex in $p$.
\eee
\end{theorem}

The next result says that Algorithm 1 terminates correctly because
any optimal solution to OPF-$\mathG$ will satisfy exactly one of the four cases
in Algorithm 1 under assumption A1.
\begin{theorem}\label{thm:algfeasible}
Suppose A1 holds.   Given an optimal solution $x^*$ of OPF-$\mathG$, exactly one of the
following holds:
\bee
\item[C1]: $P^*_{0,1}\leq 0$.
\item[C2]: $P^*_{n,0'}\geq 0$.
\item[C3]: $\exists ! k\in \mathN$ such that $P^*_{k,k+1}\geq 0$ and $P^*_{k+1,k}\geq 0$.
\item[C4]: $\exists ! k\in \mathN$ such that $P^*_{k,k-1}\leq 0$  and $P^*_{k,k+1}\leq 0$.
\eee
\end{theorem}

The intuition behind Theorem \ref{thm:algfeasible} is that if more than one of C1-C4 are satisfied, there will be at least one load bus in $(1\dots n)$ that injects positive real power, which violates A1. Theorem \ref{thm:exact} and \ref{thm:algfeasible} guarantee that Algorithm 1 is feasible and terminates correctly under A1-A4. Next, we will study the suboptimality gap of Algorithm 1. Some of the structure properties that will be used to proved the results are relegated to Appendix \ref{app:structure}.

To obtain the suboptimality bound of Algorithm 1, we need to define an OPF problem as sequel. 
\bq\label{eqn:opfs}
\text{OPF-$\mathG$s:}\quad \quad f(p_0):=\min_{x\in \mathbb{X}(\mathG)} p_{0'}\quad\mbox{s.t. } p_0 \mbox{ is a constant}
\eq

Based on Theorem \ref{thm:algfeasible}, there exists a unique solution $x^*$ for any OPF problems with convex objective function under A2 and A3. Hence there is also a unique solution to OPF-$\mathG s$ $x^*$ for any feasible real power injection $p_0$ at substation $0$. In other words, $x^*$ is a function of $p_0$ and let  $x(p_0;\mathG):=(p,q,P,Q)$ represents the solution to OPF-$\mathG s$ with  real power injection $p_0$ at substation $0$. We skip $v$ and $\ell$ in $x$ since $v_i$ is fixed by A2 and $\ell_{ij}$ is uniquely determined by $P_{ij}$ and $Q_{ij}$ according to \eqref{eqn:lastopf}. By Maximum theorem, $x(p_0;\mathG)$ is a continuous function of $p_0$.

Let $I_{p_0}:=\{p_0\mid \exists x\in\mathbb{X}(\mathG)\}$ represent the projection of $\mathbb{X}(\mathG)$ on real line. $I_{p_0}$ is compact since $\mathbb{X}(\mathG)$ is compact by A4. $f(p_0)$ is strictly convex and monotone decreasing by Corollary \ref{coro:convexf}, it is right differentiable and denote its right derivative by $f_+'(p_0)$, which is monotone increasing and right differentiable and denote its right derivative by $f^{''}_{++}(p_0)$. Let
\bq\label{eqn:curvature1}
\kappa_f:=\inf_{p_0\in I_{p_0}}f^{''}_{++}(p_0)\geq 0.
\eq
$\kappa_f$ represents the minimal value of the curvature on a compact interval if $f(p_0)$ is twice differentiable.

\begin{figure}
\centering
\subfloat[Same voltage magnitude] {
	\includegraphics[scale=0.22]{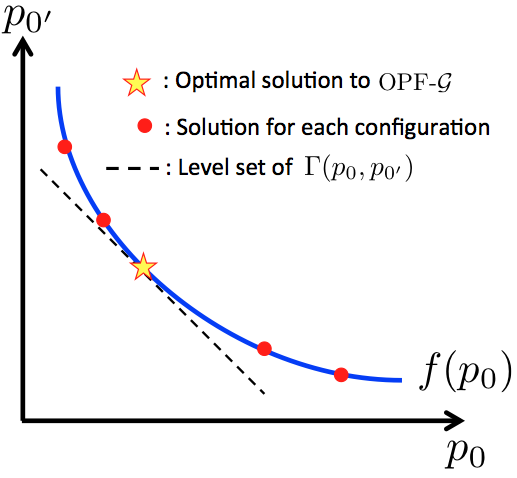}
	\label{fig:vsame}
	}
\subfloat[Different voltage magnitudes]{
	\includegraphics[scale=0.22]{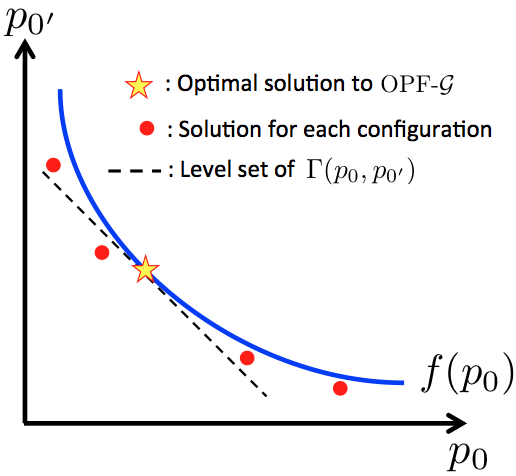}
	\label{fig:vdifferent}
	}
% \caption{Two newly connected feeders served by different substations (left panel)
%	or the same substation (right panel).}
\caption{Graph interpretations for Theorem \ref{thm:vsame} and \ref{thm:bound}}
\label{fig:intuition}
\end{figure}

Before formally state the result, we will first explain the intuition of the suboptimality gap using Fig. \ref{fig:intuition}. Since OPF-$\mathG$ can be written equivalently as 
\bqn
\min \ \ \Gamma(p_0,f(p_{0})),
\eqn
solving OPF-$\mathG$ is equivalent to find a point when the level set of $\Gamma(p_0,p_{0'})$ first hits the curve $(p_0,f(p_0))$ on a two dimensional plane, where the x-axis and y-axis are the real power injections from substation $0$ and $0'$, as shown in Fig. \ref{fig:intuition}. On the other hand, OFR can be written as 
\bqn
\min_{(i,j)\in\mathE(0,0')} \ \ \Gamma(p_0^i,p_{0'}^j)
\eqn
and solving OFR is equivalent to find a point when the level set of $\Gamma(p_0,p_{0'})$ first hits one point in $\{(p_0^i,p_{0'}^j)\mid (i,j)\in\mathE(0,0')\}$ on the two dimentional plane.

 When the voltage magnitude of each bus is fixed at the same value, all the $(p_0^i,p_{0'}^j)$ locates exactly on the curve $(p_0,f(p_0))$ as shown in Fig. \ref{fig:vsame}. Thus, we can obtain exactly the optimal solution to OFR by checking the points $(p_0^i,p_{0'}^j)$ adjacent to the optimal solution to OPF-$\mathG$. When the voltage magnitude of each bus is different, $(p_0^i,p_{0'}^j)$ does not locate on the curve $(p_0,f(p_0))$ as shown in Fig. \ref{fig:vdifferent}. Thus, the points $(p_0^i,p_{0'}^j)$ adjacent to the optimal solution to OPF-$\mathG$ may not be optimal for OFR. However, we show that the suboptimality gap is related to three aspects: 1) the distance of $(p_0^i,p_{0'}^j)$ to the curve $(p_0,f(p_0))$ (depicted by line loss), 2) the distance between each point in $(p_0^i,p_{0'}^j)$ (depicted by the power injection at each bus) and 3) the convexity of $f(p_0)$ (depicted by the curvature $\kappa_f$). Since the line loss is much smaller than  the power injection at load buses, the suboptimality gap is usually very small, as discussed after Theorem \ref{thm:bound}.
 
% However, the distance of $(p_0^i,p_{0'}^j)$ to the curve $(p_0,f(p_0))$ (depicted by line loss) is much smaller than the distance between each point in $(p_0^i,p_{0'}^j)$ (depicted by the power injection at each bus). Then, the points adjacent to the optimal solution to OPF-$\mathG$ are very close to the optimal solution to OFR since $f(p_0)$ is convex and a more convex $f(p_0)$ (larger $\kappa_f$) gives a smaller suboptimality bound.

Now, we will formally state our results on the suboptimality bound of Algorithm 1. When the voltage magnitude of all the buses are fixed at the same reference value, e.g. $1$ p.u., Algorithm 1 finds an optimal solution to OFR.
\begin{theorem}\label{thm:vsame}
Suppose A1--A4 hold.  If the voltage magnitudes of all buses are fixed at the same
value, then the line $e^*$ returned by Algorithm 1 is optimal for OFR.
\end{theorem}

When the voltage magnitudes are fixed but different at different buses,
Algorithm 1 is not guaranteed to find a global optimum of OFR.
However it still gives an excellent suboptimal solution to OFR. 
By nearly optimal, it means the suboptimality gap of Algorithm 1 is negligible.

%gap between optimal solution to OFR and solution given by Algorithm 1 is bounded above by a small value.
% Furthermore, the assumption that the reactive power is not binding is also reasonable since reactive power is usually free and can be provided by either placing inverters or capacitors on each bus.

Define $L_k$ for each line $(k,k+1)\in \mathE(0,0')$ as sequel.
\bqn
L_k:=\frac{\delta v_{k}^2 \, r_{k,k+1}/|z_{k,k+1}|^2}{(v_k+v_{k+1})+\sqrt{(v_k+v_{k+1})^2-\delta v_{k}^2
			\left( \frac{r^2_{k,k+1}}{x^2_{k,k+1}}+1 \right)}}
\eqn
where $\delta v_k:=v_k-v_{k+1}$. 
$L_k$ represents the thermal loss of line $(k,k+1)$ when either $P_{k,k+1}$ or $P_{k+1,k}$ is $0$. 
Conceptually it means all the real power sending from bus on one end of the line is converted to thermal loss 
and the other bus receives $0$ real power, namely either $P_{k,k+1}=\ell_{k,k+1}r_{k,k+1}$ or $P_{k+1,k}=\ell_{k,k+1}r_{k,k+1}$.
Then the expression of $L_k=\ell_{k,k+1}r_{k,k+1}$ can be obtained by substituting either $P_{k,k+1}=\ell_{k,k+1}r_{k,k+1}$ or $P_{k+1,k}=\ell_{k,k+1}r_{k,k+1}$ into \eqref{eqn:opf2} and \eqref{eqn:lastopf}. $L_k$ is negligible compared to the power consumption of a load in a distribution system. 
Therefore the ratio of these two quantity, defined as $R_{k}:=-\overline p_{k+1}/L_k$, is usually quite large.
\slow{Why is $L_k$ the thermal loss?   Is it equivalent to $\ell_{k,k+1}r_{k,k+1}$?  Is it easy to see?}
\qiuyu{It needs some algebra. I have put more comment as to how to get the expression above.}

Let $R:=\min R_k$ and $\kappa_f$ as defined in \eqref{eqn:curvature1}, which is a constant depending on the network. Let $\Gamma^*$ be the optimal objective value of OFR and $\Gamma_{A}$ be the objective value if we open the line $e^*$  given by Algorithm 1.
 \begin{theorem}\label{thm:bound}
Suppose A1--A4 hold and for all $i\in \mathN$.
Then
\bqn
\Gamma^*\leq \Gamma_{A}\leq \Gamma^*+\max\left\{\frac{c_0^2}{c_{0'}},\frac{c_{0'}^2}{c_{0}}\right\}\frac{2}{R^2\kappa_f},
\eqn
if $\Gamma(p_0,p_{0'}):=c_0p_0+c_{0'}p_{0'}$ for some positive $c_0,c_{0'}$.
\end{theorem}

{\it Remark}: $R$ is large, usually on the order of $10^3$, in a distribution system when there is no renewable generation. Although it is difficult to estimate the value of $\kappa_f$ in theory, our simulation shows that $\kappa_f$ is typically around $0.025 MW^{-1}$ for a feeder with loop size of $10$, thus the bound is approximately $80W$ if $c_0=c_{0'}=1$, which is quite small. Moreover simulations of two SCE distribution circuits show that Algorithm 1  always finds the global optima of OFR problem; see section \ref{sec:simulation}. Therefore the bound in the theorem, already negligible, is not always tight.

The suboptimality bound depends on the expression of the objective function $\Gamma(p_0,p_{0'})$ if it is not linear but strictly convex. According to Theorem \ref{thm:bound}, we know that a more convex $f(p_0)$ (larger $\kappa_f$) gives a smaller suboptimality gap. In general, a more convex objective function also suggests a smaller suboptimality gap, which can also be interpreted using Fig. \ref{fig:vdifferent}.

\qiuyu{The algorithm solves OFR if $|p_0^{i+1}-p_0^i|>\frac{1}{R\kappa_f}$ if $c_0=c_{0'}$. $\frac{1}{R\kappa_f}\approx 0.04MW$ in our simulation and that explains why our algorithm returns the optimal solution.}

%------------------------------------------------------------------------------------------------------------------------------------------------------------------------

\section{Simulation}\label{sec:simulation}
\begin{figure*}
\centering
\includegraphics[width=0.9\linewidth]{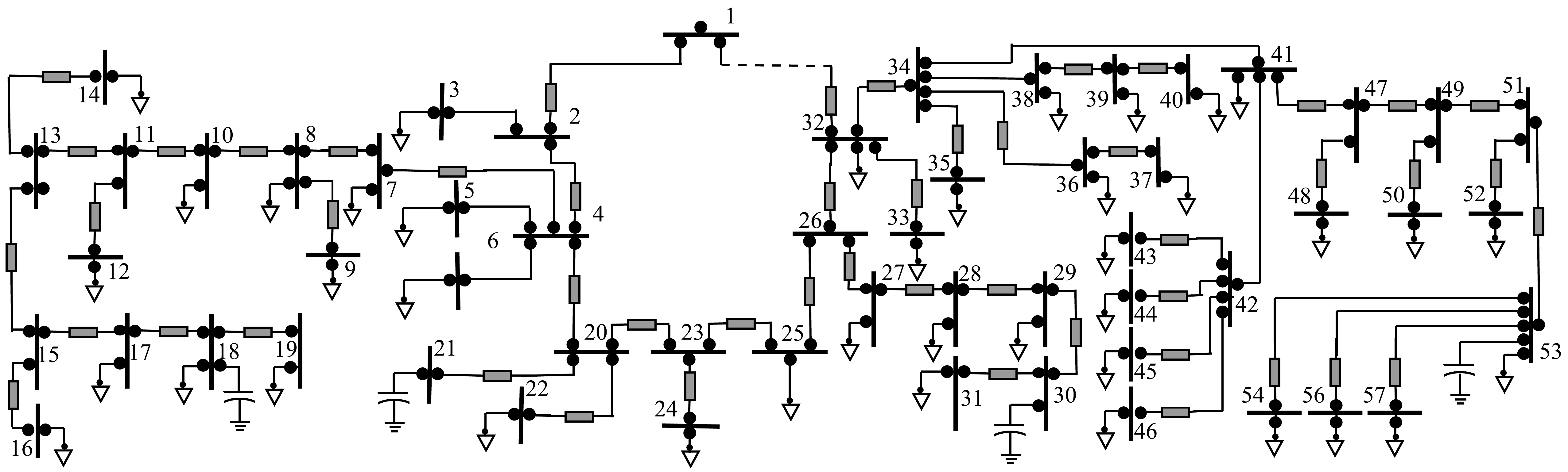}
\caption{A modified SCE $56$ bus feeder. Bus $1$ is the substation and line $(1,32)$ is added.}
\label{fig:56bus}
\end{figure*}

In this section we present an example to illustrate the effectiveness of Algorithm 1. 
The simulation is implemented using the CVX optimization toolbox \cite{cvx2012} in Matlab. We use a $56$-bus 
SCE distribution feeder whose circuit diagram is shown in Fig. \ref{fig:56bus}. 
The network data, including line
 impedances and real power demand of loads, are listed in Table \ref{table:data}. 
 Since there is no loop in the  original feeder we added a tie line between bus 1 and bus 32, which is
 assumed to be initially open.
 
 \begin{table}
\caption{The aggregate power injection from substation $1$ for each configuration}
\centering
\begin{tabular}{|c|c|c|c|c|}
\hline
Opened line & $(1,2)$ & $(2,4)$  & $(4,20)$ & $(20,23)$\\
\hline
Power injection (MW) & $3.8857$ & $3.8845$  & $3.8719$ & $3.8718$\\
\hline
Opened line & $(23,25)$ & $(25,26)$  & $(26,32)$ & $(32,1)$\\
\hline
Power injection (MW) & $3.8719$ & $3.8721$  & $3.8755$ & $3.9550$\\
\hline
\end{tabular}
\label{tab:case1}
\end{table}

\begin{figure}
\centering
\includegraphics[width=0.9\linewidth]{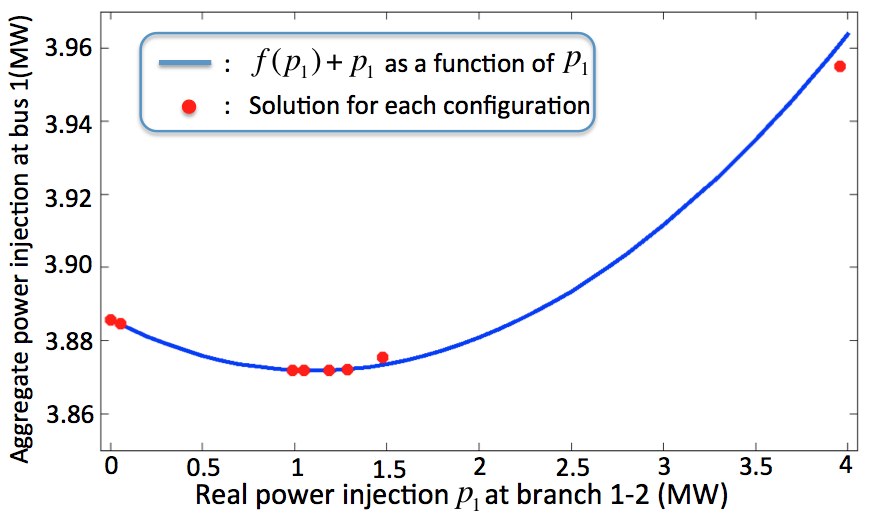}
\caption{The aggregated power injection $f(p_1)+p_1$ as a function of $p_1$. The red dots are the 
operating points for each configuration.}
\label{fig:case1}
\end{figure}

In our simulation the voltage magnitude of the substation (bus $1$) is fixed at 1 p.u.. 
We relax the assumption A2 needed for our analysis that their voltage magnitudes at all other buses are fixed
and allow them to vary within $[0.97,1.03]$p.u., as required in the current distribution system. The demand of real power is fixed for each load and shown in Table \ref{table:data}. The reactive power at each bus, which is kept within $10\%$ of the real power to maintain a 
power factor of at least $90\%$, is a control variable, as in Volt/VAR control of \cite{Farivar-2011-VAR-SGC}.

\begin{table*}
\caption{Line impedances, load demand and capacitors for the distribution circuit in Figure \ref{fig:56bus}. }
\centering
\scriptsize
\begin{tabular}{|c|c|c|c|c|c|c|c|c|c|c|c|c|c|c|c|c|c|}
\hline
\multicolumn{18}{|c|}{Network Data}\\
\hline
\multicolumn{4}{|c|}{Line Data}& \multicolumn{4}{|c|}{Line Data}& \multicolumn{4}{|c|}{Line Data}& \multicolumn{2}{|c|}{Load Data}& \multicolumn{2}{|c|}{Load Data}&\multicolumn{2}{|c|}{Load Data}\\
\hline
From&To&R&X&From&To& R& X& From& To& R& X& Bus& Peak & Bus& Peak &  Bus &Peak  \\
Bus.&Bus.&$(\Omega)$& $(\Omega)$ & Bus. & Bus. & $(\Omega)$ & $(\Omega)$ & Bus.& Bus.& $(\Omega)$ & $(\Omega)$ & No.&  MW& No.& MW& No.& MW\\
\hline

1	&	2	&	0.160	&	0.388	&	20	&	21	&	0.251	&	0.096	&	39	&	40	&	2.349	&	 0.964	&	3	&	0.057	&	29  &	0.044  &   50   & 0.045 \\
2	&	3	&	0.824	&	0.315	&	21	&	22	&	1.818	&	0.695	&	34	&	41	&	0.115	&	 0.278	&	5	&	0.121	&	31	&	0.053  &    52& 0.315	 \\
2	&	4	&	0.144	&	0.349	&	20	&	23	&	0.225	&	0.542	&	41	&	42	&	0.159	&	 0.384	&	6	&	0.049	&	32	&	0.223 &   54& 	 0.061	 \\
4	&	5	&	1.026	&	0.421	&	23	&	24	&	0.127	&	0.028	&	42	&	43	&	0.934	&	 0.383	&	7	&	0.053	&	33	&	0.123 & 	 55&	 0.055	 \\
4	&	6	&	0.741	&	0.466   &	23	&	25	&	0.284	&	0.687	&	42	&	44	&	0.506	&	 0.163	&	8	&	0.047	&	34	&	0.067 &    56&	 0.130  \\\cline{17-18}
4	&	7	&	0.528	&	0.468	&	25	&	26	&	0.171	&	0.414	&	42	&	45	&	0.095	&	 0.195	&	9	&	0.068	&	35	&	0.094&     \multicolumn{2}{c|}{Shunt Cap}	 \\\cline{17-18}
7	&	8	&	0.358	&	0.314	&	26	&	27	&	0.414	&	0.386	&	42	&	46	&	1.915	&	 0.769	&	10	&	0.048	&	36	&	0.097&     \multicolumn{1}{c|}{Bus} &	 \multicolumn{1}{c|}{Mvar}			 \\\cline{17-18}
8	&	9	&	2.032	&	0.798	&	27	&	28	&	0.210	&	0.196	&	41	&	47	&	0.157	&	 0.379	&	11	&	0.067	&	37	&	0.281&    19& 	 0.6 	 \\
8	&	10	&	0.502	&	0.441	&	28	&	29	&	0.395	&	0.369	&	47	&	48	&	1.641	&	 0.670	&	12	&	0.094	&	38	&	0.117&     21&	0.6 	 \\
10	&	11	&	0.372	&	0.327	&	29	&	30	&	0.248	&	0.232	&	47	&	49	&	0.081	&	 0.196	&	14	&	0.057	&	39	&	0.131&     30&	      0.6 		 \\
11	&	12	&	1.431	&	0.999	&	30	&	31	&	0.279	&	0.260	&	49	&	50	&	1.727	&	 0.709	&	16	&	0.053	&	40	&	0.030& 	53    &	0.6 		 \\\cline{17-18}
11	&	13	&	0.429	&	0.377	&	26	&	32	&	0.205	&	0.495	&	49	&	51	&	0.112	&	 0.270	&	17	&	0.057	&	41	&	0.046&         \multicolumn{2}{c|}{} \\
13	&	14	&	0.671	&	0.257	&	32	&	33	&	0.263	&	0.073	&	51	&	52	&	0.674	&	 0.275	&	18	&	0.112	&	42	&	0.054&     \multicolumn{2}{c|}{} \\
13	&	15	&	0.457	&	0.401	&	32	&	34	&	0.071	&	0.171	&	51	&	53	&	0.070	&	 0.170	&	19	&	0.087	&	43	&	0.083&    \multicolumn{2}{c|}{$V_\textrm{base}$ = 12kV}  \\
15	&	16	&	1.008	&	0.385	&	34	&	35	&	0.625	&	0.273	&	53	&	54	&	2.041	&	 0.780	&	22	&	0.063	&	44	&	0.057&  \multicolumn{2}{c|}{$S_\textrm{base}$ = 1MVA} \\
15	&	17	&	0.153	&	0.134	&	34	&	36	&	0.510	&	0.209	&	53	&	55	&	0.813	&	 0.334	&	24	&	0.135	&	45	&	0& 	\multicolumn{2}{c|}{$Z_\textrm{base}= 144 \Omega$ } \\
17	&	18	&	0.971	&	0.722	&	36	&	37	&	2.018	&	0.829	&	53	&	56	&	0.141	&	 0.340	&	25	&	0.100	&	46	&	0.134&  	\multicolumn{2}{c|}{} \\
18	&	19	&	1.885	&	0.721	&	34	&	38	&	1.062	&	0.406	&	32	&	1	&	0.085  	&	 0.278	 &	27	&	0.048	&	47	&	0.045&  	\multicolumn{2}{c|}{} \\
4	&	20	&	0.138	&	0.334	&	38	&	39	&	0.610	&	0.238	&		&		&		&		 &	28	&	0.038	&	48	&	0.196&  		\multicolumn{2}{c|}{} 	 \\
		
\hline
\end{tabular}
\label{table:data}
\end{table*}

We use the aggregate power injection as our objective, namely $\Gamma(p_0,p_{0'}):=p_0+p_{0'}$.
It also represents the power loss in this case since we have fixed real power demand of each load. 
Our addition of the line between buses $1$ and $32$ creates a loop 1-2-4-20-23-25-26-32-1 that
must be broken by turning off the switch on one line from $\{(1,2),(2,4),(4,20),(20,23),(23,25),$ $(25,26),(26,32),(32,1)\}$. 
In Table \ref{tab:case1} we list the corresponding aggregate power injection for all the possible configurations.
The optimal configuration is to open line $(20,23)$ at an optimal cost of $3.8718$ MW.

After we run Algorithm 1 bus $23$ receives real power from both sides and our algorithm returns line $(20,23)$, 
which is the optimal solution to OFR. 

We have explained the suboptimality gap derived in Theorem \ref{thm:bound} using Fig. \ref{fig:vdifferent}. We claim that the suboptimality gap is small because the distance of the point $(p_0^i,p_{0'}^j)$ to the curve $(p_0,f(p_0))$ is small for real system, which suggests a negligible gap. Now we verify the claim using our simulation results. Let $p_1$ and $f(p_1)$ denote the real power flow on branch $(1,2)$ and 
$(32,1)$ respectively. Instead of plotting $f(p_1)$ versus $p_1$ as Fig. \ref{fig:vdifferent}, we plot the aggregate power injection $f(p_1)+p_1$ 
at substation $1$ as a function of $p_1$ for $p_1\in[0,4]$ in Fig. \ref{fig:case1} for better illustration\footnote{We actually run the topology where 
the substation is virtually broken into two substations as Fig. \ref{fig:step1a}.}. In Fig. \ref{fig:case1}, the red dots corresponding to the operational points for each configuration in Table \ref{tab:case1} and the blue curve is $f(p_1)+p_1$. All the red dots are much closer to the curve than their adjacent points, which illustrates why the suboptimality gap is small in real system. 

\slow{What is the point you'd like to illustrate?  Perhaps be more explicit?}

Even though the voltage magnitudes in our simulation are not fixed at the nominal values as 
assumed in our analysis, Algorithm 1 still gives the optimal solution to OFR by solving a convex
relaxation of OPF. The underlying reason is that the voltage magnitude does not vary much 
between adjacent buses in real network, hence the performance of the algorithm is not limited by 
the assumption of fixed voltage magnitudes. 

\slow{The logic here seems weak?   Even if voltages across each line differ by a small amount,
the difference can accumulate and it is possible to have very different voltages at the head and
at the end of a feeder. We should either explain better or remove the explanation and simply 
point out the fact.}
We have also tested our algorithm in another SCE $47$ bus distribution feeder and it
again yields the optimal solution to OFR.

%-----------------------------------------------------------------------------------------------------------------------------------------------------------------------------------------------

%-----------------------------------------------------------------------------------------------------------------------------------------------------------------------------------------------
\section{Conclusion}\label{sec:conclusion}

We have proposed an efficient algorithm to optimize the branch exchange step in 
feeder reconfiguration, based on SOCP relaxation of OPF.
We have derived a bound on the suboptimality gap and argued that it is very small.
We have proved that the algorithm computes an optimal solution when all voltage magnitudes
are the same.  
We have demonstrated the effectiveness of our algorithm through simulations of real-world 
feeders.

%----------------------------------------------------------------------------------------------------------------------------------------------------------------------------------------------
%\input{tex/8_bio}
 \bibliographystyle{IEEEtran}
 \bibliography{IEEEabrv,tex/reference,tex/PowerRef-201202}

\begin{thebibliography}{10}
\providecommand{\url}[1]{#1}
\csname url@samestyle\endcsname
\providecommand{\newblock}{\relax}
\providecommand{\bibinfo}[2]{#2}
\providecommand{\BIBentrySTDinterwordspacing}{\spaceskip=0pt\relax}
\providecommand{\BIBentryALTinterwordstretchfactor}{4}
\providecommand{\BIBentryALTinterwordspacing}{\spaceskip=\fontdimen2\font plus
\BIBentryALTinterwordstretchfactor\fontdimen3\font minus
  \fontdimen4\font\relax}
\providecommand{\BIBforeignlanguage}[2]{{%
\expandafter\ifx\csname l@#1\endcsname\relax
\typeout{** WARNING: IEEEtran.bst: No hyphenation pattern has been}%
\typeout{** loaded for the language `#1'. Using the pattern for}%
\typeout{** the default language instead.}%
\else
\language=\csname l@#1\endcsname
\fi
#2}}
\providecommand{\BIBdecl}{\relax}
\BIBdecl

\bibitem{Merlin1975}
A.~Merlin and H.~Back, ``Search for a minimal-loss operating spanning tree
  configuration in an urban power distribution system,'' in \emph{Proc. of the
  Fifth Power System Conference (PSCC), Cambridge}, 1975, pp. 1--18.

\bibitem{Civanlar1988}
S.~Civanlar, J.~Grainger, H.~Yin, and S.~Lee, ``Distribution feeder
  reconfiguration for loss reduction,'' \emph{Power Delivery, IEEE Transactions
  on}, vol.~3, no.~3, pp. 1217--1223, 1988.

\bibitem{Baran1989c}
M.~Baran and F.~Wu, ``{Network reconfiguration in distribution systems for loss
  reduction and load balancing},'' \emph{IEEE Trans. on Power Delivery},
  vol.~4, no.~2, pp. 1401--1407, Apr 1989.

\bibitem{Chiang1990a}
H.-D. Chiang and R.~Jean-Jumeau, ``Optimal network reconfigurations in
  distribution systems: Part 1: A new formulation and a solution methodology,''
  \emph{IEEE Trans. Power Delivery}, vol.~5, no.~4, pp. 1902--1909, November
  1990.

\bibitem{sarfi1994}
R.~J. Sarfi, M.~Salama, and A.~Chikhani, ``A survey of the state of the art in
  distribution system reconfiguration for system loss reduction,''
  \emph{Electric Power Systems Research}, vol.~31, no.~1, pp. 61--70, 1994.

\bibitem{Jabr2011}
R.~El~Ramli, M.~Awad, and R.~Jabr, ``Ordinal optimization for dynamic network
  reconfiguration,'' \emph{Electric Power Components and Systems}, vol.~39,
  no.~16, pp. 1845--1857, 2011.

\bibitem{Cormen2001}
T.~H. Cormen, C.~E. Leiserson, R.~L. Rivest, and C.~Stein, \emph{Introduction
  to Algorithms}, 2nd~ed.\hskip 1em plus 0.5em minus 0.4em\relax The MIT Press,
  2001.

\bibitem{Morton2000}
A.~B. Morton and I.~M. Mareels, ``An efficient brute-force solution to the
  network reconfiguration problem,'' \emph{Power Delivery, IEEE Transactions
  on}, vol.~15, no.~3, pp. 996--1000, 2000.

\bibitem{Liu1989}
C.-C. Liu, S.~J. Lee, and K.~Vu, ``Loss minimization of distribution feeders:
  optimality and algorithms,'' \emph{IEEE Trans. Power Delivery}, vol.~4,
  no.~2, pp. 1281--1289, April 1989.

\bibitem{Dall'Anese2013-reconfig1}
E.~Dall'Anese and G.~B. Giannakis, ``Sparsity-leveraging reconfiguration of
  smart distribution systems,'' August 2013, arXiv:1303.5802v2.

\bibitem{Dall'Anese2013-reconfig2}
------, ``Risk-constrained microgrid reconfiguration using group sparsity,''
  June 2013, arXiv:1306.1820v1.

\bibitem{Wang1995}
J.-C. Wang, H.-D. Chiang, and G.~R. Darling, ``An efficient algorithm for
  real-time network reconfiguration in large scale unbalanced distribution
  systems,'' in \emph{Power Industry Computer Application Conference, 1995.
  Conference Proceedings., 1995 IEEE}.\hskip 1em plus 0.5em minus 0.4em\relax
  IEEE, 1995, pp. 510--516.

\bibitem{branch1}
M.~E. Baran and F.~F. Wu, ``Optimal capacitor placement on radial distribution
  systems,'' \emph{Power Delivery, IEEE Transactions on}, vol.~4, no.~1, pp.
  725--734, 1989.

\bibitem{branch2}
M.~Baran and F.~F. Wu, ``Optimal sizing of capacitors placed on a radial
  distribution system,'' \emph{Power Delivery, IEEE Transactions on}, vol.~4,
  no.~1, pp. 735--743, 1989.

\bibitem{Farivar-2011-VAR-SGC}
M.~Farivar, C.~R. Clarke, S.~H. Low, and K.~M. Chandy, ``Inverter var control
  for distribution systems with renewables,'' in \emph{Smart Grid
  Communications (SmartGridComm), 2011 IEEE International Conference on}.\hskip
  1em plus 0.5em minus 0.4em\relax IEEE, 2011, pp. 457--462.

\bibitem{Farivar-2013-BFM-TPS}
M.~Farivar and S.~H. Low, ``Branch flow model: relaxations and convexification
  (parts {I, II}),'' \emph{IEEE Trans. on Power Systems}, vol.~28, no.~3, pp.
  2554--2572, August 2013.

\bibitem{Gan-2013-BFMt-CDC}
L.~Gan, N.~Li, U.~Topcu, and S.~H. Low, ``Optimal power flow in distribution
  networks,'' in \emph{Proc. 52nd IEEE Conference on Decision and Control},
  December 2013, in arXiv:12084076.

\bibitem{V_opf}
J.~Lavaei, D.~Tse, and B.~Zhang, ``Geometry of power flows and optimization in
  distribution networks,'' \emph{arXiv preprint arXiv:1204.4419}, 2012.

\bibitem{cvx2012}
M.~Grant and S.~Boyd, ``{CVX}: Matlab software for disciplined convex
  programming, version 2.0 beta,'' \url{http://cvxr.com/cvx}, September 2012.

\bibitem{Bose-2012-BFMe-Allerton}
S.~Bose, S.~H. Low, and M.~Chandy, ``Equivalence of branch flow and bus
  injection models,'' in \emph{50th Annual Allerton Conference on
  Communication, Control, and Computing}, October 2012.

\bibitem{convex}
S.~P. Boyd and L.~Vandenberghe, \emph{Convex optimization}.\hskip 1em plus
  0.5em minus 0.4em\relax Cambridge university press, 2004.

\end{thebibliography}

%-----------------------------------------------------------------------------------------------------------------------------------------------------------------------------------------------
\appendices

\section{Proof of Theorem \ref{thm:exact}}
\noindent {\it Proof of Part 1:} 
SOCP relaxation for bus injection model has been shown in \cite{V_opf} to be exact under A2 and A3.
Since the bus injection model and the branch flow model are equivalent \cite{Bose-2012-BFMe-Allerton},
this implies that SOPF-$\mathG$ is exact under A2 and A3.

\noindent{\it Proof of Part 2:} We will prove it by contradiction. 
By part 1), SOPF-$\mathG$ is exact.  Hence every solution of SOPF-$\mathG$ is also feasible for 
OPF-$\mathG$. Assume there are two optimal solutions $\hat x$ and $\overline x$ such that they achieve the same objective value, namely $\Gamma({\overline p})=\Gamma({\hat p})$. In addition, $v_i\hat\ell_{ij}=\hat P_{ij}^2+\hat Q_{ij}^2$ and $v_i\overline\ell_{ij}=\overline P_{ij}^2+\overline Q_{ij}^2$ for $(i,j)\in \mathE$ since the relaxation is exact. Let $x(\lambda):=\lambda \hat x+(1-\lambda) \overline x$ for $0\leq\lambda\leq 1$. Then $x(\lambda)$ is feasible for SOPF-$\mathG$ since $\mathbb{X}_c(\mathG)$ is convex. Since $\Gamma(p)$ is convex,
\bqn
\Gamma(p(\lambda))\leq \Gamma({\overline p})=\Gamma({\hat p}).
\eqn
 Note that $\Gamma({\overline p})$ and $\Gamma({\hat p})$ are optimal for OPF-$\mathG$ by assumption, the only possibility is $\Gamma(p(\lambda))=\Gamma({\overline p})=\Gamma({\hat p})$, which means $x(\lambda)$ is also an optimal solution. However 
\bqn
&&P_{ij}^2(\lambda)+Q_{ij}^2(\lambda)\\
&=&(\lambda \hat P_{ij}+(1-\lambda) \overline P_{ij})^2+(\lambda \hat Q_{ij}+(1-\lambda) \overline Q_{ij})^2\\
&\leq& \lambda(\hat P_{ij}^2+\hat Q_{ij}^2)+(1-\lambda)(\overline P_{ij}^2+\overline Q_{ij}^2)\\
&=&v_i\ell_{ij}(\lambda)
\eqn
The equality is attained for any $0<\lambda< 1$ if and only if $\hat P_{ij}=\overline P_{ij}$ and $\hat Q_{ij}=\overline Q_{ij}$. Since the convex relaxation is exact, $P_{ij}^2(\lambda)+Q_{ij}^2(\lambda)=v_i\ell_{ij}(\lambda)$, which indicates $\hat\ell_{ij}=\overline \ell_{ij}$. Finally, we have $\hat s_i=\overline s_i$ by  \eqref{eqn:opf3} and \eqref{eqn:opf4}. Therefore $\hat x=\overline x$, which contradicts that there are more than one optimal solutions.

\label{app:exact}

%-----------------------------------------------------------------------------------------------------------------------------------------------------------
\section{Proof of Theorem \ref{thm:algfeasible}}

For any solution $x^*$ to OPF-$\mathG$, define
\bqn
\mathcal{I}_P^*:=\{P^*_{0,1},-P^*_{1,0},P^*_{1,2},-P^*_{2,1},\ldots, P^*_{n,0'},-P^*_{0',n}\}. 
\eqn
We will show that the elements in $\mathcal{I}_P^*$ are in a descending order and
at most two consecutive elements are equal.   This implies that exactly one of C1--C4 holds.

By power balance across each line $(i,i+1)\in \mathE(0,0')$ we have
\bqn
P^*_{i,i+1}=-P^*_{i+1,i}+\ell^*_{i,i+1}r_{i,i+1}\geq -P^*_{i+1,i}.
\eqn
By power balance at each bus $i\in \mathN\setminus\{0,0'\}$ we have
\bqn
P^*_{i,i+1}=-P^*_{i,i-1}+p^*_i<-P^*_{i,i-1}
\eqn
under assumption A1 because $p^*_i\leq \overline p_i<0$. Hence, $\mathcal{I}^*_{P}$ is a nonincreasing sequence and we are left to show that there exists at most two equal element in $\mathcal{I}^*_{P}$. We will show it by contradiction. Since $P^*_{i,i+1}$ is strictly less than $-P^*_{i,i-1}$ under A1, the equality can only happens between $P^*_{i,i+1}$ and $-P^*_{i+1,i}$. Suppose there exists two lines $(k_1,k_1+1),(k_2,k_2+1)\in \mathE(0,0')$ such that $P^*_{k_1,k_1+1}=-P^*_{k_1+1,k_1}$ and $P^*_{k_2,k_2+1}=-P^*_{k_2+1,k_2}$. It means $\ell_{k_1,k_1+1}=\ell_{k_2,k_2+1}=0$, which indicates $P^*_{k_1+1,k_1}=P^*_{k_2,k_2+1}=0$ by \eqref{eqn:lastopf}. Assume $k_1<k_2$ without loss of generality and by power balance equation \eqref{eqn:opf3}, 
\bqn
P^*_{k_1+1,k_1}+P^*_{k_2,k_2+1}=\sum_{i=k_1+1}^{k_2}p_i^*-\sum_{i=k_1+1}^{k_2-1}\ell^*_{i,i+1}r_{i,i+1}<0,
\eqn
which contradicts $P^*_{k_1+1,k_1}=P^*_{k_2,k_2+1}= 0$. Thus, there are at most two equal elements.

%-----------------------------------------------------------------------------------------------------------------------------------------------------------------------------------------------
\section{Proof of Theorem  \ref{thm:vsame}}\label{app:vsame}

By Lemma \ref{lem:concave}, $P_{k,k+1}(p_0;\mathG)$ is an increasing and continuous function of $p_0$, hence there exists a unique $p_0$ to $P_{k,k+1}(p_0;\mathG)=0$ and let it to be $p_0({k})$. Recall that $p_0^k$ and $p_{0'}^{k+1}$ defined in \eqref{eqn:opfi} and \eqref{eqn:opfj}, we will first show $p_0^k=p_0(k)$ and $p_{0'}^{k+1}=f(p_0(k))$ for any $(k,k+1)\in \mathcal{E}(0,0')$ provided the voltage magnitude of each bus is the same. By symmetry, It suffices to show $p_0^k=p_0(k)$. $P_{k,k+1}(p_0(k);\mathG)=0$ indicates $Q_{k,k+1}(p_0(k);\mathG)=0$ according to \eqref{eqn:opf2}. Hence, $p_0(k)$ is a feasible power injection for subtree $\mathG_0^k$ and it means $p_0^k\leq p_0(k)$. Next, we will show that it is the smallest possible power injection for $\mathG_0^k$. Suppose we have $p_0^k<p_0(k)$, then $(p_0^k,f(p_0(k)))$ is a feasible power injection for network $\mathG$ with $p_0^k<p_0(k)$. It contradicts $(p_0(k),f(p_0(k)))\in\mathcal{O}(\mathbb{P})$. Therefore we have $p_0^k=p_0(k)$, $p_{0'}^{k+1}=f(p_0(k))$ and
\bqn
(p_0^k,p_{0'}^{k+1})=(p_0(k),f(p_0(k)))\in\mathcal{O}(\mathbb{P}),
\eqn
which means the minimal power injection for each partition of graph $\mathG$ locates exactly on the Pareto front of the feasible power injection region of OPF-$\mathG$. Therefore OFR is equivalent to the following problem:
\bqn
\min_{0\leq k\leq n}\Gamma(p_0^k,p_{0'}^{k+1})
=\min_{0\leq k\leq n}\Gamma(p_0(k),f(p_0(k))),
\eqn
whose minimizer is denoted by $k^*$. Similarly, OPF-$\mathG$ can be rewritten as
\bqn
\min_{p_0\in I_{p_0}} \Gamma(p_0,f(p_0)),
\eqn
whose unique minimizer is denoted by $p_0^*$ and let $x^*$ be the optimal solution. By Lemma \ref{lem:concave}, $P_{k,k+1}(p_0;\mathG)$ is an increasing function of $p_0$. Therefore we have
\bee
\item $P_{0,1}(p_0^*;\mathG)<0$ $\Leftrightarrow$ $k^*=0$.
\item $P_{0',n}(p_0^*;\mathG)<0$ $\Leftrightarrow$ $k^*=n$.
\item $\exists k\in[0,n]$ such that $-P_{k,k-1}(p_0^*;G)\leq 0$ and $P_{k,k+1}(p_0^*;\mathG)\geq 0$ $\Leftrightarrow$ $k^*=k-1$ or $k$.
\eee
The above $3$ cases correspond to 1), 2) and 4) in Algorithm 1. Note that 3) would never happen when all the bus voltages are fixed at the same magnitude.

%-----------------------------------------------------------------------------------------------------------------------------------------------------------------------------------------------
\section{Proof of Theorem \ref{thm:bound}}\label{sec:theorem4}

Let $p_0^f({k})$ be the solution to $P_{k,k+1}(p_0;\mathG)=0$ and $p_0^b({k})$ be the solution to $P_{k+1,k}(p_0;\mathG)=0$. The uniqueness of $p_0^f({k})$ and $p_0^b({k})$ can be shown in a similar manner as the uniqueness of $p_0(k)$ in the proof of Theorem \ref{thm:vsame}.
When the voltage magnitude of each bus is the same, $p_0^f({k})=p_0^b({k})$ and they degenerate to $p_0(k)$.
\begin{lemma}\label{lem:qfree}
Suppose A1-A4 hold. For any $(k,k+1)\in \mathcal{E}(0,0')$, $p_0^k=p^f_0(k)$ and $p_{0'}^{k+1}=f(p^b_0(k))$.
\end{lemma}

\begin{IEEEproof}
It suffices to show $p_0^f(k)$ is optimal for $\mathG_0^k$ due to symmetry. First, we show $p_0^f(k)$ is feasible for $\mathG_0^k$. Given a solution  $x(p_0^f(k);\mathG)$ to OPF-$\mathG s$, let $\tilde x:=(\tilde P,\tilde Q,\tilde p,\tilde q)$, where
\bqn
\tilde S_{i,i+1}&=&S_{i,i+1}(p_0^f(k);\mathG) \quad i< k\\
%\tilde \ell_{i,i+1}&=&\ell_{i,i+1}(p_0^f(k);G) \quad i< k\\
\tilde s_{i}&=&s_i(p_0^f(k);\mathG)\quad i<k\\
\tilde p_k&=&p_k(p_0^f(k);\mathG) \\
 \tilde q_k&=&q_k(p_0^f(k);\mathG)-Q_{k,k+1}(p_0^f(k);\mathG)
\eqn
Thus, we have $\tilde x\in \mathbb{X}(\mathG_0^k)$, which means $p_0^f(k)$ is feasible for $\mathG_0^k$. Next, we will show $p_0^f(k)$ is the minimal power injection for $\mathG_0^k$. Suppose $\hat p_0<p_0^f(k)$ is feasible for $\mathG_0^k$, then we can construct a feasible solution $\tilde x\in \mathbb{X}(\mathG)$ and the real power injection at node $0$ and $0'$ are $\hat p_0$ and $f(p_0^f(k))$, respectively. Therefore it contradicts that $(p_0^f(k),f(p_0^f(k)))\in\mathcal{O}(\mathbb{P})$. The construction process is as follows: 
\bqn
\tilde S_{i,i+1}=
\begin{cases}
S_{i,i+1}(\hat p_0;\mathG) & i< k\\
S_{i,i+1}(p_0^f(k);\mathG) & i\geq k
\end{cases}
\eqn
\bqn
\tilde s_{i}=
\begin{cases}
s_i(\hat p_0;\mathG) & i< k\\
p_k(\hat p_0;\mathG)+ \mathbf{i}(q_k(\hat p_0;\mathG)+Q_{k,k+1}(p_0^f(k);\mathG)) & i=k\\
s_i(p_0^f(k);\mathG) & i> k
\end{cases}
\eqn
It can be verified that $\tilde x \in \mathbb{X}(\mathG)$. Therefore, we show that $p_0^f(k)$ is the minimal power injection for $\mathG_0^k$.
\end{IEEEproof}

\begin{lemma}\label{lem:bound}
Suppose A1-A4 hold. Then we have
\bqn
\frac{p_0^b(k)-p_0^f(k)}{p_0^f(k+1)-p_0^b(k)}\leq \frac{1}{R_{k}}.
\eqn
\end{lemma}
\begin{IEEEproof}
By Lemma \ref{lem:concave}, $P_{k,k+1}(p_0,\mathG)$ is a concave increasing function with respect to $p_0$. Hence $-P_{k+1,k}(p_0,\mathG)=\phi(P_{k,k+1}(p_0,\mathG))$ is also a concave increasing function of $p_0$. Recall that $-P_{k+1,k}(p_0^f(k),\mathG)=-L_k$, $-P_{k+1,k}(p_0^b(k),\mathG)=0$ and $-P_{k+1,k}(p_0^f(k+1),\mathG)=-p_{k+1}\geq -\overline p_{k+1}$, we have
\bqn
\frac{0-(-L_{k})}{p_0^b(k)-p_0^f(k)}\geq\frac{-\overline p_{k+1}-0}{p_0^f(k+1)-p_0^b(k)}
\eqn
by definition of a concave function. Rearrange the above inequality, we obtain
\bqn
\frac{p_0^b(k)-p_0^f(k)}{p_0^f(k+1)-p_0^b(k)}\leq \frac{0-(-L_{k})}{-\overline p_{k+1}-0}:=\frac{1}{R_k}.
\eqn
\end{IEEEproof}

\begin{lemma}\label{lem:math1}
Let $g(x)$ be a strictly convex decreasing function supported on $[a,b]$ and $\kappa_g:=\inf_{x\in(a,b)}g^{''}_{++}(p_0)$. Define $G(x):=c_1g(x)+c_2x$ $(c_1,c_2>0)$, which is also strictly convex with a unique minimizer $x^*$ on $[a,b]$. Let $a\leq y_1\leq \cdots\leq y_{2n-1}\leq y_{2n}\leq b$ be a partition on $[a,b]$ such that
\bq\label{eqn:101}
\frac{y_{2i}-y_{2i-1}}{y_{2i+1}-y_{2i}}\leq \frac{1}{R} \qquad (1\leq i\leq n-1)
\eq
for some $R>0$. Then there exists a $0\leq k\leq 2n$ such that $y_{k}\leq x^*\leq y_{k+1}$, where $y_{0}=a$ and $y_{2n+1}=b$. Let $G_i:=c_1g(y_{2i})+c_2y_{2i-1}$ for $1\leq i\leq n$ and $G^*:=\min_{1\leq i\leq n}\left\{G_i\right\}$.
Define
\begin{displaymath}
G_A:=
\begin{cases}
G_1& \text{if } k=0\\
G_n &\text{if } k=2n\\
G_{(k-1)/2} &\text{if } k\text{ is odd}\\
\min\{G_{k/2},G_{k/2+1}\} &\text{if } k\neq 0,2n \text{ and is even}
\end{cases}
\end{displaymath}
Then
\bqn
G^*\leq G_{A}\leq G^*+\max\left\{\frac{c_1^2}{c_{2}},\frac{c_{2}^2}{c_{1}}\right\}\frac{2}{R^2\kappa_g}.
\eqn
\end{lemma}
\begin{IEEEproof}
Without loss of generality, assume $c_1\leq c_2$ and let $\lambda:=\frac{c_2}{c_1}$. The unique minimizer $x^*$ of $G(x):=c_1g(x)+c_2x$ is
\bqn
x^*&=&\arg\sup_{x\in[a,b]}\{x\mid G'_+(x)\leq 0\}\\
&=&\arg\sup_{x\in[a,b]}\{x\mid g'_+(x)\leq -\lambda\}.
\eqn
In addition, let
\bqn
x_l&:=&\arg\sup_{x\in[a,b]}\left\{x\mid g'_+(x)\leq-\lambda(1+\frac{1}{R})\right\}\\
x_r&:=&\arg\inf_{x\in[a,b]}\left\{x\mid g'_+(x)\geq-\lambda(1-\frac{1}{R})\right\}
\eqn
Then we have  $x_l\leq x^*\leq x_r$ because $g(x)$ is strictly convex. Let
\begin{align*}
t_1&:=\max\{i\mid y_{2i-1}\leq x^*\}  \quad t_2:=\min\{i\mid y_{2i}\geq 2x_l-x^*\}\\
t_3&:=\min\{i\mid y_{2i}\geq x^*\}  \quad t_4:=\max\{i\mid y_{2i-1}\leq 2x_r-x^*\}
\end{align*}

Next, we will prove the result for different $k$ as sequel.

\noindent {Case I: $k=2n$.} In this case, we have $d_a=G_n$, $t_1=n$. We need to further divide it into two categories.

(1.a) [$y_{2t_1-1}\leq x_l$ or $y_{2t_1-1}\in[x_l,x^*]$ and $t_1=t_2$]. \\
For any $i<t_1$,
\bq
 \frac{G_i-G_{i+1}}{c_0}
&=&\lambda g(y_{2i})+y_{2i-1}-g(y_{2i+2})-\lambda y_{2i+1}\nonumber\\
&\geq& \lambda g(y_{2i})+y_{2i-1}-g(y_{2i+1})-\lambda y_{2i+1}\nonumber\\
&=& \lambda (y_{2i-1}-y_{2i})+G(y_{2i})-G(y_{2i+1})\nonumber\\
&\geq& \frac{\lambda}{R}(y_{2i}-y_{2i+1})+G(y_{2i})-G(y_{2i+1})\nonumber\\
&\geq& 0 \label{eqn:102}
\eq
The first inequality follows from $g(x)$ is an increasing function and the second inequality follows from the assumption \eqref{eqn:101}. For the last inequality, if $y_{2t_1-1}\leq x_l$, we have $G'_+(y_{2i+1})<-\lambda/R$ for all $i<n$ and the inequality holds according to mean value theorem. If $y_{2t_1-1}\in[x_l,x^*]$ and $t_1=t_2=n$, the inequality holds for $i<n-1$ due to similar reason above. When $i=n-1$,
\bqn
&&G(y_{2i})-G(y_{2i+1})\\
&\geq& G'_+(x_l)(y_{2i}-y_{2i+1})\geq -\frac{\lambda}{R}(y_{2i}-y_{2i+1})
\eqn
because $y_{2n-2}\leq 2x_l-x^*$ by definition of $t_2$ and $G(x)$ is convex. \eqref{eqn:102} means the sequence $\{G_i\}$ is of descending order and $G^*=G_n=G_A$, thus $G^*-G_A=0$.

(1.b) [$y_{2t_1-1}\in[x_l,x^*]$ and $t_2<t_1$].\\
 In this case, $y_{2i+1}-y_{2i}\leq 2(x^*-x_l)$ for $t_2\leq i< t_1$. Denote
\bq\label{eqn:103}
\delta y_{i}:=y_{2i}-y_{2i-1}\leq\frac{y_{2i+1}-y_{2i}}{R}\leq \frac{2(x^*-x_l)}{R}
\eq
 for $t_2\leq i< t_1$. Note that the curvature of $g(x)$ is bounded below by $\kappa_g$, then $x^*-x_l\leq \lambda/(R\kappa_g)$. Substitute it into \eqref{eqn:103}, we have for $t_2\leq i< t_1$
 \bq\label{eqn:104}
 \delta y_i\leq \frac{2\lambda}{R^2\kappa_g}.
 \eq
Then for $t_2\leq i< t_1$,
\bqn
G_i-G_{t_1}&=& c_1g(y_{2i})+c_2y_{2i-1}-G_{t_1}\\
&\geq &c_1g(y_{2i})+c_2y_{2i-1}-G(y_{2t_1-1})\\
&=& -c_2\delta y_i+G(y_{2i})-G(y_{2t_1-1})\\
&\geq & -\frac{2c_2^2}{c_1R^2\kappa_g}\\
\eqn
Clearly the first inequality holds. The second inequality follows from  \eqref{eqn:104} and $G(x)$ is monotone decreasing for $x\leq x^*$.

For $i\leq t_2$, $G_i>G_{t_1}$ can be shown in a similar manner as (1.a). Thus, we have $G_i-G_{n}\geq -\frac{2c_2^2}{c_1R^2\kappa_g}$ for any $i\leq n$, which indicates $G_A-G^*\leq -\frac{2c_2^2}{c_1R^2\kappa_g}$.

\noindent {Case II: $k=0$.} In this case, $G_A=G_1$ and the bound can be established in a similar manner as Case I.

\noindent {Case III: $k$ is odd.} In this case, $G_A=G_{(k-1)/2}$, $t_1=t_3=(k-1)/2$. Similar approach can be applied as Case I and Case II to show
\bqn
G_i\geq
\begin{cases}
G(y_{2t_1-1})-\frac{2c_2^2}{c_1R^2\kappa_g} &\text{ if } i\leq t_1\\
G(y_{2t_3})-\frac{2c_2^2}{c_1R^2\kappa_g}&\text{ if } i\geq t_3
\end{cases}
\eqn
And $G_A=G_{(k-1)/2}\leq \max\{G(y_{2t_1-1}),G(y_{2t_3})\}\leq G_i+\frac{2c_2^2}{c_1R^2\kappa_g}$ for $1\leq i\leq n$.

\noindent{Case IV: $k\neq 0$ and is even.} In this case, $G_A=\min\{G_{k/2},G_{k/2+1}\}$ and $t_1=k/2$, $t_3=k/2+1$. Similar approach can be applied as Case I and Case II to show
\bqn
G_i\geq
\begin{cases}
G_{t_1}-\frac{2c_2^2}{c_1R^2\kappa_g} &\text{ if } i\leq t_1\\
G_{t_3}-\frac{2c_2^2}{c_1R^2\kappa_g} &\text{ if } i\geq t_3
\end{cases}
\eqn
and we arrive at our conclusion.
\end{IEEEproof}

Consider the sequence $p_0^f(0)\leq p_0^b(0)\leq \ldots\leq p_0^f(n)\leq p_0^b(n)$ as the partition on $I_{p_0}$  and $f(p_0)$ as the function $g(x)$ in Lemma \ref{lem:math1}, we can prove Theorem \ref{thm:bound}.

%---------------------------------------------------------------------------------------------------------------------------------------------------------------------------------------------------

\section{Structural Properties of OPF}\label{app:structure}

Given two real vectors $x,y\in\mathbb{R}^n$, $x\leq y$ means $x_i\leq y_i$ for $1\leq i\leq n$ and $x<y$ means $x_i<y_i$ for at least one component. The Pareto front of a compact set $A\subseteq \mathbb{R}^n$ is defined as
\bqn
\mathcal{O}(A):=\{x\in A\mid \nexists \tilde x\in A\setminus\{x\} \text{ such that } \tilde x\leq x\}
\eqn
Let $\mathbb{P}:=\{(p_0,p_{0'})\mid \exists x\in\mathbb{X}(\mathG)\}$ represent the projection of $\mathbb{X}(\mathG)$ on $\mathbb{R}^2$ and $\mathbb{P}_c:=\{(p_0,p_{0'})\mid \exists x\in\mathbb{X}_c(\mathG)\}$ be the projection of $\mathbb{X}_c(\mathG)$ on $\mathbb{R}^2$. Based on Theorem \ref{thm:exact}, the convex relaxation is exact under A2 and A3, it means  $\mathcal{O}(\mathbb{P})=\mathcal{O}(\mathbb{P}_c)$ since the objective $\Gamma(p_0,p_{0'})$ is a convex and nondecreasing function. Next, we begin by studying some properties of $f(p_0)$ defined in \eqref{eqn:opfs}.

\begin{lemma}\label{lem:equiv1}
$(p_0,f(p_0))\in \mathcal{O}(\mathbb{P})$ provided the optimization problem OPF-$\mathG s$ is feasible.
\end{lemma}

\begin{IEEEproof}
Suppose $(p_0,\hat{p}_{0'})\in \mathcal{O}(\mathbb{P})$. Let $\Gamma^*(p_0,p_{0'})$ be the objective function such that $(p_0,\hat{p}_{0'})$ solves OPF-$\mathG$. Note that  OPF-$\mathG$ can be written equivalently as
\bqn
\min_{p_0} \Gamma^*(p_0,p_{0'}) \quad \text{s.t. }p_{0'}=f(p_0)
\eqn
Therefore, at optimality, $f(p_0)=\hat{p}_{0'}$ and $(p_0,f(p_0))\in \mathcal{O}(\mathbb{P})$.
\end{IEEEproof}

By property of {\it Pareto Front} \cite{convex}, for each point $(p_0,p_{0'})\in\mathcal{O}(\mathbb{P}_c)=\mathcal{O}(\mathbb{P})$, there exists a convex nondecreasing function $\Gamma^*:\mathbb{R}^2\rightarrow \mathbb{R}$ such that $(p_0,p_{0'})$ is an optima for OPF-$\mathG$. Therefore OPF-$\mathG s$ and OPF-$\mathG$ are equivalent in the sense that if we fix $p_0$ and solve OPF-$\mathG s$ with $x^*$, there exists an objective function $\Gamma^*(p_0,p_{0'})$ such that  $x^*$ solves the corresponding OPF-$\mathG$ by Lemma \ref{lem:equiv1}.

\begin{lemma}\label{lem:cd}
Let $A$ be a compact and convex set in $\mathbb{R}^2$. Define $g(x):=y$ for  any $(x,y)\in \mathcal{O}({A})$. Then $y=g(x)$ is a convex decreasing function of $x$ for $(x,y)\in \mathcal{O}({A})$.
\end{lemma}

\begin{IEEEproof}
We first show $g(x)$ is a decreasing function and then show $g(x)$ is also convex.

Let $(x_1,g(x_1))$ and $(x_2,g(x_2))$ be two points in $\mathcal{O}({A})$. Without loss of generality, assume $x_1>x_2$. If $g(x_1)\geq g(x_2)$, it violates the fact that $(x_1,g(x_1))\in \mathcal{O}({A})$ and hence $g(x_1)< g(x_2)$, which means that $g(x)$ is a decreasing function.

Next, we will show $g(\cdot)$ is convex. Recall that $A$ is a compact set, we have $(x_1,g(x_1)),(x_2,g(x_2))\in \mathcal{O}({A})\subseteq A$. $A$ is also a convex set, thus $(\frac{x_1+x_2}{2},\frac{g(x_1)+g(x_2)}{2})\in A$. By definition of Pareto front,
\bqn
g(\frac{x_1+x_2}{2})=\inf_{(\frac{x_1+x_2}{2},y)\in A}\left\{y\right\}\leq  \frac{g(x_1)+g(x_2)}{2},
\eqn
which shows $g(x)$ is a convex function.
\end{IEEEproof}

Note that $\mathbb{X}_c(\mathG)$ is convex and compact by A4, hence its projection on a two dimensional space $\mathbb{P}_c$ is also compact and convex. By Lemma \ref{lem:equiv1}, $(p_0,f(p_0))$ characterizes the Pareto front $\mathcal{O}(\mathbb{P})=\mathcal{O}(\mathbb{P}_c)$, then we have the following corollary.
\begin{corollary}\label{coro:convexf}
$f({p_0})$ is a strictly convex decreasing function of $p_0$ under assumption A2-A4.
\end{corollary}

For a line $(i,j)$ between two buses $i$ and $j$ with fixed voltage magnitude, $(P_{ij},Q_{ij},\ell_{ij})$ are governed by \eqref{eqn:opf2}-\eqref{eqn:lastopf} and $Q_{ij},\ell_{ij}$ can be uniquely solved given a $P_{ij}$ if A3 holds. Let $\phi(P_{ij}):=-P_{ji}=P_{ij}-\ell_{ij}r_{ij}$ and we have the following result.

\begin{lemma}\label{lem:branchincrease}
Suppose A2 and A3 hold, $\phi(P_{ij})$ is a concave increasing function of $P_{ij}$ for $(i,j)\in \mathE$.
\end{lemma}

\begin{IEEEproof}
By \eqref{eqn:lastopf}, we have $\ell_{ij}=(P_{ij}^2+Q_{ij}^2)/v_i$ and substitute it in $\phi(P_{ij})$, we have
\bqn
\phi(P_{ij})=P_{ij}-\frac{r_{ij}}{v_i}\left(P_{ij}^2+Q_{ij}^2\right).
\eqn
The relation between $P_{ij}$ and $Q_{ij}$ is governed by \eqref{eqn:opf2}. Let $\theta_{ij}:=\theta_i-\theta_j$ and then
$P_{ij}$ and $Q_{ij}$ can be written as
\bqn
P_{ij}&=&\frac{v_ir_{ij}}{r_{ij}^2+x_{ij}^2}+\sqrt{\frac{v_iv_j}{r_{ij}^2+x_{ij}^2}}\sin(\theta_{ij}-\beta_{ij})\\
Q_{ij}&=&\frac{v_ix_{ij}}{r_{ij}^2+x_{ij}^2}-\sqrt{\frac{v_iv_j}{r_{ij}^2+x_{ij}^2}}\cos(\theta_{ij}-\beta_{ij}),
\eqn
where $\beta_{ij}:=\arctan r_{ij}/x_{ij}$. Substitute them into $\phi(P_{ij})$, we obtain
\bqn
\phi(P_{ij})=-\frac{v_jr_{ij}}{r_{ij}^2+x_{ij}^2}+\sqrt{\frac{v_iv_j}{r_{ij}^2+x_{ij}^2}}\sin(\theta_{ij}+\beta_{ij}).
\eqn
Take derivative of $\phi(P_{ij})$ with respect to $P_{ij}$, we have
\bqn
\frac{d\phi(P_{ij})}{d P_{ij}}=\frac{\cos(\theta_{ij}+\beta_{ij})}{\cos(\theta_{ij}-\beta_{ij})}.
\eqn
which is always positive by assumption A3 that $|\theta_{ij}|< \arctan x_{ij}/r_{ij}$. Furthermore,
\bqn
%\frac{d^2\phi(P_{ij})}{d P_{ij}^2}=\frac{2r_{ij}x_{ij}}{\sqrt{(r_{ij}^2+x_{ij}^2)v_iv_j}\sin^3\theta_{ij}}<0,
\frac{d^2\phi(P_{ij})}{d P_{ij}^2}=-\sqrt{\frac{r_{ij}^2+x_{ij}^2}{v_iv_j}}\frac{\sin2\beta_{ij}}{\cos^3(\theta_{ij}-\beta_{ij})},
\eqn
which is always negative by assumption A3 that $|\theta_{ij}|< \arctan x_{ij}/r_{ij}$. Thus, $\phi(P_{ij})$ is a concave increasing function of $P_{ij}$.

\end{IEEEproof}

Lemma \ref{lem:branchincrease} means if the one end of the line increases its real power injection on the line, the other end should receive more real power under assumption A2 and A3.

%\begin{figure}
%\centering
%\includegraphics[width=0.9\linewidth]{figure/lemma45}
%\caption{Illustration of Lemma \ref{lem:increase} and \ref{lem:concave}}
%\label{fig:lemma45}
%\end{figure}

\begin{lemma}\label{lem:increase}
Suppose A2-A4 hold. Given a solution $x(p_0;\mathG)$ to OPF-$\mathG s$, $P_{k,k+1}(p_0;\mathG)$ is a nondecreasing function of $p_0$ for all $(k,k+1)\in \mathE(0,0')$.
\end{lemma}

\begin{IEEEproof}
The following argument holds without assuming $\overline q_i=\infty$ and $\underline q_i=-\infty$.

Suppose $P_{k,k+1}(p_0;\mathG)$ is not a nondecreasing function of $p_0$ at $p_0^*$ for a line $(k,k+1)\in\mathE(0,0')$, then either C1 or C2 below will hold for arbitrary small $\epsilon>0$,

\begin{itemize}
\item [C1:] $\exists p_0\in(p_0,p_0^*+\epsilon)$ such that $P_{k,k+1}(p_0;\mathG)<P_{k,k+1}(p_0^*;\mathG)$.
\item [C2:] $\exists p_0\in(p_0-\epsilon,p_0^*)$ such that $P_{k,k+1}(p_0;\mathG)>P_{k,k+1}(p_0^*;\mathG)$.
\end{itemize}

We will show by contradiction that $(p^*_0,f(p^*_0))\not \in \mathcal{O}(\mathbb{P})$ in this case, which violates Lemma \ref{lem:equiv1}. Assume without loss of generality that $P_{i,i+1}(p_0,\mathG)$ is a nondecreasing function of $p_0$ for $0\leq i<k$.

\noindent{\it Case I: $q_k(p^*_0;\mathG)> \underline q_k$}. Suppose C1 holds, then there exists a monotone decreasing sequence $p_0^{(m)}\downarrow p^*_0$ such that 
$\{P_{k,k+1}(p_0^{(m)};\mathG),m\in\mathbb{N}\}$ is a monotone increasing sequence that converges to $P_{k,k+1}(p_0^*;\mathG)$ because $x(p_0;\mathG)$ is continuous over $p_0$. By power balance equation \eqref{eqn:opf3} at bus $k$, for any $m$, we have
\bqn
p_k(p_0^{(m)};\mathG)&=&P_{k,k+1}(p_0^{(m)},\mathG)-\phi(P_{k-1,k}(p_0^{(m)};\mathG))\\
&<&P_{k,k+1}(p^{(m+1)}_0,\mathG)-\phi (P_{k-1,k}( p^{(m+1)}_0,\mathG))\\
&=&p_k(p^{(m+1)}_0;\mathG)
\eqn
Thus $\{p_k(p_0^{(m)};\mathG),n\in\mathbb{N}\}$ is a monotone increasing sequence that converges to $p_k(p_0^*;\mathG)$. We now construct a point $\tilde x=(\tilde P,\tilde Q,\tilde p,\tilde q)$ as follows. First, pick up $(\tilde P_{k,k+1},\tilde Q_{k,k+1},\tilde p_k,\tilde q_k)$ such that  $\tilde p_k\in\{p_k(p_0^{(m)};\mathG ),m\in\mathbb{N}\}$, $\tilde q_k\in( \underline q_k, q_k(p^*_0;\mathG))$ and they satisfy the following equations:
\begin{subequations}\label{eqn:app300}
\begin{align}
 \tilde P_{k,k+1}=&P_{k,k+1}(p^*_0;\mathG)-p_k(p^*_0;\mathG)+\tilde p_k\\
 \tilde Q_{k,k+1}=&Q_{k,k+1}(p^*_0;\mathG)-q_k(p^*_0;\mathG)+\tilde q_k\\
 v_{k+1}=&v_k-2(r_{k,k+1}\tilde P_{k,k+1}+x_{k,k+1}\tilde Q_{k,k+1})\\
 &+\frac{\tilde P_{k,k+1}^2+\tilde Q_{k,k+1}^2}{v_k}|z_{k,k+1}|^2
\end{align}
\end{subequations}
The existence of $(\tilde P_{k,k+1},\tilde Q_{k,k+1},\tilde p_k,\tilde q_k)$ is guaranteed by the following two facts:
\bi
\item $p_k(p_0^{(m)};\mathG)$ is a monotone increasing sequence that converges to $p_k(p_0^*;\mathG)$.
\item $\tilde q_k$ is a continuous decreasing function of $\tilde p_k$ if they satisfy \eqref{eqn:app300}.
\ei
Since $\tilde P_{k,k+1}\in[P_{k,k+1}(p^{(1)}_0;\mathG),P_{k,k+1}(p^*_0;\mathG)]$ and $x(p_0;\mathG)$ are continuous over $p_0$, then there exists a $ p'_0\in[p^*_0,p^{(1)}_0]$ such that $S_{k,k+1}(p'_0;\mathG)=\tilde S_{k,k+1}$.

Next, we will construct the feasible physical variable for $i\neq k$. For $0\leq i<k$, let $\tilde s_i=s_i(p^*_0;\mathG)$ and $\tilde S_{i,i+1}=S_{i,i+1}(p^*_0;\mathG)$. For $k<i\leq n$, let $\tilde s_i=s_i(p'_0;\mathG)$ and $\tilde S_{i,i+1}=S_{i,i+1}(p'_0;\mathG)$. Clearly that $\tilde x\in\mathbb{X}(\mathG)$ with $(p^*_0,f(p'_0))$ as the real power injection at substation $0$ and $0'$. However, $f(p'_0)<f(p^*_0)$, which contradicts $(p^*_0,f(p^*_0))\in \mathcal{O}(\mathbb{P})$.

\noindent{\it Case II: $q_k(p^*_0;G)<\overline q_k$}. Similar approach can be used to show C2 does not hold by contradiction.

So far, we have shown that $P_{k,k+1}(p_0;\mathG)$ is non-decreasing either on its left or right neighborhood. Thus $P_{k,k+1}(p_0;\mathG)$ is non-decreasing of $p_0$ if $\underline q_k<\overline q_k$ because $P(p_0;\mathG)$ is a continuous function of $p_0$. The case where $\underline q_k=\overline q_k$ can be covered by taking limitation of the case of $\underline q_k<\overline q_k$.

\end{IEEEproof}

\begin{lemma}\label{lem:concave}
Suppose A2-A4 hold. Given a solution $x(p_0;\mathG)$ to OPF-$\mathG s$, $P_{k,k+1}(p_0;\mathG)$ is a concave increasing function of $p_0$ for all $(k,k+1)\in \mathE(0,0')$.
\end{lemma}

\begin{IEEEproof}
It is shown that $P_{k,k+1}(p_0;\mathG)$ is a nondecreasing function of $p_0$ in Lemma \ref{lem:increase}. We now show it is also a concave function of $p_0$. Let $\mathG_1=(\mathN_1,\mathE_1)$, where $\mathN_1=\{i\mid 0\leq i\leq k+1\}$ and $\mathE_1=\{(i,i+1)\mid 1\leq i\leq k\}$. All the physical constraints are the same as $\mathG$ except the bus injection power $s_{k+1}$ at node $k+1$, which is relaxed to be a free variable. Mathematically, it means
\bqn
\overline \ell_{i,i+1}(\mathG_1)&=&\overline \ell_{i,i+1}(\mathG)\quad i\leq k \\
\overline s_{i}(\mathG_1)&=&\overline s_i(\mathG) \quad i\leq k  \qquad \overline s_{k+1}(\mathG_1)=\infty \\
\underline s_i(\mathG_1)&=&\underline s_i(\mathG)\quad i\leq k  \qquad \underline s_{k+1}(\mathG_1)=-\infty
\eqn
Consider the following OPF problem:
\bqn
\text{OPF-$\mathG s1$:}\qquad \quad \min_{x\in\mathbb{X}(\mathG_1)}p_{k+1} \quad\mbox{s.t. } p_0 \mbox{ is a fixed}
\eqn
Let $p_{k+1}(p_0;\mathG_1)$ be the optimal value for OPF-$\mathG s1$ and $P_{k,k+1}(p_0;\mathG_1)$ be the real branch power flow across line $(k,k+1)$, respectively.

Next, we will show $P_{k,k+1}(p^*_0;\mathG)=P_{k,k+1}(p^*_0;\mathG_1)$. Clearly that $P_{k,k+1}(p^*_0;\mathG)\leq P_{k,k+1}(p^*_0;\mathG_1)$. Otherwise, by Lemma \ref{lem:branchincrease}, we have
\bqn
-p_{k+1}(p_0^*;\mathG_1)=\phi(P_{k,k+1}(p^*_0;\mathG_1))< \phi(P_{k,k+1}(p^*_0;\mathG)),
\eqn
which contradicts that $p_{k+1}(p_0^*;\mathG_1)$ is optimal for OPF-$\mathG s1$. Thus, it suffices to show $P_{k,k+1}(p^*_0;\mathG)< P_{k,k+1}(p^*_0;\mathG_1)$ does not hold. Suppose $P_{k,k+1}(p^*_0;\mathG)< P_{k,k+1}(p^*_0;\mathG_1)$ holds. By Lemma \ref{lem:increase}, $P_{k,k+1}(p_0;\mathG)$ is a nondecreasing function of $p_0$, thus there exists a $\hat p_0>p_0^*$ such that $P_{k,k+1}(\hat p_0;\mathG)\in[P_{k,k+1}(p^*_0;\mathG),P_{k,k+1}(p^*_0;\mathG_1))$. Recall that $\mathbb{X}_c(\mathG)$, which is the set of feasible solutions after the SOCP relaxation, is convex and is connected, there exists a $x\in\mathbb{X}_c(\mathG_1)$ with $P_{k,k+1}=P_{k,k+1}(\hat p_0;\mathG)$ but $p_0=p_0^*$. It means $(p^*_0,f(\hat p_0))$ is also feasible for OPF-$\mathG$, which contradicts that $(\hat p_0,f(\hat p_0))\in\mathcal{O}(\mathbb{P})$.

Now we have $P_{k,k+1}(p^*_0;\mathG)= P_{k,k+1}(p^*_0;\mathG_1)$. Since the convex relaxation is exact, $p_{k+1}(p_0;\mathG_1)$ is a convex decreasing function of $p_0$ by Lemma \ref{lem:cd}. In addition,
\bqn
\phi(P_{k,k+1}(p_0;\mathG))=\phi(P_{k,k+1}(p_0;\mathG_1))=-p_{k+1}(p_0,\mathG_1),
\eqn
where $\phi(\cdot)$ is a continuous increasing function, thus is invertible. Then $P_{k,k+1}(p_0,\mathG)=\phi^{-1}(-p_{k+1}(p_0,\mathG_1))$ is a concave function of $p_0$.
\end{IEEEproof}

\end{document}